\documentclass[11pt]{article}

\usepackage[draft]{ed}
\usepackage{amsfonts}
\usepackage{amsmath}
\usepackage[]{algorithm2e}
\usepackage{amssymb}
\usepackage{epsfig}
\usepackage{bm}
\usepackage{enumerate}
\numberwithin{equation}{section} 

\usepackage{color}
\usepackage[normalem]{ulem}

\title{Multitrace formulations and Domain Decomposition Methods for the solution of Helmholtz transmission problems for bounded composite scatterers}
\author{Carlos Jerez-Hanckes, Carlos P\'erez-Arancibia, Catalin Turc}

\newtheorem{theorem}{Theorem}[section]

\newtheorem{remark}[theorem]{Remark}
\newenvironment{proof}{\hspace{0.5cm} {\bf Proof.}}
{$\quad {}_\blacksquare$\vspace{0.3cm}}

\date{}
\newcommand{\triple}[1]{{\left\vert\kern-0.25ex\left\vert\kern-0.25ex\left\vert #1 
    \right\vert\kern-0.25ex\right\vert\kern-0.25ex\right\vert}}

\begin{document}
\maketitle
\begin{abstract}
  We present Nystr\"om discretizations of multitrace formulations and non-overlapping Domain Decomposition Methods (DDM) for the solution of Helmholtz transmission problems for bounded composite scatterers with piecewise constant material properties. We investigate the performance of DDM with both classical Robin and generalized Robin boundary conditions. The generalized Robin boundary conditions incorporate square root Fourier multiplier approximations of Dirichlet to Neumann operators. While the classical version of DDM is not particularly well suited for Krylov subspace iterative solvers, we show that the associated DDM linear system can be efficiently solved by hierarchical elimination via Schur complements of the Robin data. We show through numerical examples that the latter version of DDM gives rise to small numbers of Krylov subspace iterations that depend mildly on the frequency and number of subdomains. 
   
\end{abstract}

\section{Introduction}
\label{intro}

The phenomenon of electromagnetic wave scattering by bounded penetrable objects composed of several subdomains with different but constant electric permittivities is relevant for numerous applications in antenna design, diffraction gratings, and photovoltaic cells, to name but a few. It is typical in all these application areas that multiple media meet at a single point, a scenario that is referred to as \emph{triple or multiple junctions}. Standard numerical methods perform poorly when dealing with wave scattering by composite objects with piecewise constant material properties portraying large frequency ranges, as they need to resolve wave interactions with high-contrast sharp interfaces. Volumetric discretizations of these problems result in very large linear systems of equations that are ill-conditioned in the high-frequency regime and whose solution by iterative solvers require inordinate numbers of iterations.  Several preconditioning strategies have been proposed to mitigate the aforementioned issues, the most successful arguably being those that rely on the shifted Laplacean~\cite{bayliss1983iterative,erlangga2004class} or the sweeping preconditioner introduced in~\cite{engquist2011sweeping}.

Domain Decomposition Methods (DDM) are natural candidates for the solution of scattering problems involving composite scatterers. DDM are divide and conquer strategies whereby the computational domain is divided into smaller subdomains for which solutions are matched via transmission conditions on subdomain interfaces. The convergence of DDM for time-harmonic wave scattering applications depends a great deal on the choice of the transmission conditions that allow the exchange of information between adjacent subdomains. These interface transmission conditions should ideally allow information to flow out of a subdomain with as little as possible information being reflected back into the subdomain. Thus, the interface transmission conditions fall into the category of Absorbing Boundary Conditions (ABC). From this perspective, the ideal choice of transmission conditions on an interface between two subdomains is such that the impedance/transmission operator is the restriction to the common interface of the Dirichlet-to-Neumann (DtN) operator corresponding to the adjacent subdomain. Traditionally, the interface transmission conditions were chosen as the classical first order ABC outgoing Robin/impedance boundary conditions~\cite{Depres,Collino1}. The convergence of DDM with classical Robin interface boundary conditions is slow and is adversely affected by the number of subdomains. Fortunately, the convergence of DDM can be considerably improved through incorporation of ABC that constitute higher order approximations of DtN operators in the form of second order approximations with optimized tangential derivative coefficients~\cite{Gander1}, square root approximations~\cite{boubendirDDM}, or other types of non-local transmission conditions~\cite{Collino1,steinbach2011stable}. Alternatively, so-called Perfectly Matched Layers can be used at subdomain interfaces~\cite{stolk2013rapidly}. Although the use of more sophisticated ABC, as those recounted above, accelerates a great deal the convergence of DDM, the number of iterations required for convergence still grows --albeit not drastically-- with the frequency and number of subdomains. This is not entirely surprising since higher order ABC  only involve local exchange of information between adjacent subdomains, and affect to a lesser degree the global exchange of information between distant subdomains. Recent efforts have been devoted to construct ``double sweep''-type preconditioners that address the latter issue~\cite{vion2014double,zepeda2016method}. The resulting preconditioned DDM scale favorably with frequency and number of subdomains, but appear to be somewhat less effective for wave propagation problems in composite media that exhibit sharp high-contrast interfaces.   

Boundary integral equation based solvers for scattering by composite objects with piecewise constant material properties require significantly fewer unknowns than volumetric solvers as only the interfaces of material discontinuity need be discretized. The formulation of these problems in terms of robust boundary integral equations has recently received  significant interest in the community, the main achievement being the introduction of Multi-Trace Formulations (MTF)~\cite{jerez-hanckes1,jerez-hanckes2}. The derivation of one of the MTF --dubbed \emph{local}-- consists of the following steps: (1) use of Green's identities in each subdomain, whose boundary is a union of interfaces of material discontinuity, to represent the fields in that subdomain via layer potential; (2) application of Dirichlet and Neumann traces associated to that subdomain to the Green's identities; followed by (3) enforcement of the continuity conditions across interfaces to replace the identity terms in the previous steps by Dirichlet and Neumann traces of solutions in adjacent subdomains. This procedure leads to a boundary integral equation of the first kind whose unknowns are both interior and exterior Dirichlet and Neumann traces of fields on each interface and which involves, the standard Boundary Integral Operators (BIO) on each subdomain corresponding to the wavenumber associated with that subdomain. Galerkin discretizations of the local MTF using $h$-refinement order bases for relative low frequencies preconditioned by diagonal or Calder\'on preconditioners yield small Krylov subspace iteration numbers. However, as contrast ratios increase, the Calder\'on preconditioners are less effective when used for low order Galerkin discretizations. A remedy is to resort to spectral or high-order Galerkin discretizations~\cite{JPT15,JPT16}. The MTF of the second kind---or \emph{global} version---can be derived if the fields are sought in terms of suitable linear combinations of layer potentials defined on the union of all interfaces--typically referred to as the \emph{skeleton}--of material discontinuity~\cite{jerez-hanckes2,greengard1}.

In this paper, we illustrate the performance of Nystr\"om solvers based multitrace formulations and DDM solvers in the case of high-frequency scattering problems from composite high-contrast scatterers. A major advantage of MTF is the ease with which they can be incorporated into existing boundary integral equation solvers. We present in this work a straightforward extension of the Helmholtz transmission Nystr\"om solvers introduced in~\cite{dominguez2015well} to MTF. We investigate two types of simple preconditioners for the MTF of the first kind: the Calder\'on preconditioner proposed in~\cite{jerez-hanckes1} and a new preconditioner that combines a Schur complement approach with a Calder\'on preconditioner. Both preconditioners are shown to be effective for high-frequency high-contrast scattering problems from composite scatterers not depending on the number of subdomains. However, the numbers of iterations required by Nystr\"om discretizations of MTF grows considerably with the frequency and/or the contrast between subdomains, even after resorting to preconditioning. We show that in the aforementioned frequency regime DDM based on boundary integral equations can be advantageous alternatives to MTF. We investigate both DDM based on the exchange of classical Robin data between subdomain as well as DDM that incorporate general Robin/impedance boundary operators that are square root approximations of DtN maps~\cite{HJP13}. In both cases, we solve the subdomain Helmholtz equations with Robin or generalized Robin boundary conditions using well-conditioned boundary integral formulations solved by Nystr\"om discretizations. Provided the size of the subdomains is small enough--in terms of wavelengths across--, the latter problems can be solved by direct linear algebra methods. 

The numerical results presented in this paper corroborate the by now well-known fact that DDM based on classical Robin boundary conditions do not perform well as  either the frequency or/and the number of subdomains are increased. We present two possible approaches to overcome these shortcomings. The first approach consist of hierarchical elimination from the DDM linear system via Schur complements of the Robin data corresponding to interior subdomain interfaces. This idea was presented recently as a hierarchical merging of Robin-to-Robin (RtR) maps~\cite{Gillman1} in the context of scattering from variable media with smooth index of refraction, and was shown to be equivalent to hierarchical elimination of interior Robin data from the DDM linear system corresponding to multiple scattering~\cite{pedneault2016schur}. The advantage of the Schur complement DDM is that at each stage inverses of relatively small matrices need be computed, and these can be performed in a hierarchical fashion that optimizes the computational cost without compromising the inherent parallelism of DDM. The second approach we pursue is to incorporate square root approximations--using complex wavenumbers--of the DtN operators as generalized Robin operators in DDM. The operators were used successfully for regularizing boundary integral equations at high frequencies~\cite{Antoine,AntoineX,turc1,turc2,dominguez2015well,boubendir2014well} and also incorporated successfully in DDM~\cite{boubendirDDM,vion2014double}. Since for composite scatterers, adjacent domain decomposition subdomains can have different wavenumbers, we blend the approximations of DtN operators corresponding to different subdomains--with different wavenumbers as well--via smooth cutoff functions.  

Important advantages of using Nystr\"om discretizations are (a) the square root DtN approximations which are naturally defined in the Fourier domain can be implemented with ease as Fourier multipliers using FFTs and (b) the regularization strategy presented in~\cite{turc2016well} is applicable to general Robin boundary conditions that involve the square root operators. In contrast, boundary element or other type of discretizations (e.g. finite difference or finite element) require either local approximations or rational Pad\'e approximations of the square root operators, and these approximations additional require tangential derivatives~\cite{boubendirDDM,vion2014double}. In addition, neither the presence of cross points at multiple junctions nor the incorporation of square root DtN approximations as generalized Robin boundary conditions affect the accuracy of the DDM solvers. We note that the accuracy of finite difference/finite element DDM based on generalized Robin boundary conditions has not been reported in the literature~\cite{boubendirDDM,vion2014double}. The Nystr\"om DDM with generalized Robin boundary conditions requires much fewer iterations than the classical Robin boundary conditions DDM while the computational cost of the local subdomain solves for the former method  is virtually the same as that for the latter.

The structure of this paper is as follows. Section~\ref{MS1} describes the problems of scattering considered in this paper. Section~\ref{MS1} reviews both first- and second-kind MTF approaches. The proposed DDM approach is then introduced and analyzed in Section~\ref{MS2}. Section~\ref{Nystrom}, on the other hand, discusses high-order Nytstr\"om discretizations for the various integral equation formulations considered, while a variety of numerical results are shown in Section~\ref{num}. Finally,  the conclusions of this work are presented in Section~\ref{conclu}.

\parskip 2pt plus2pt minus1pt

\section{Electromagnetic transmission problems in composite domains\label{MS1}}

\begin{figure}
\centering
\includegraphics[height=60mm]{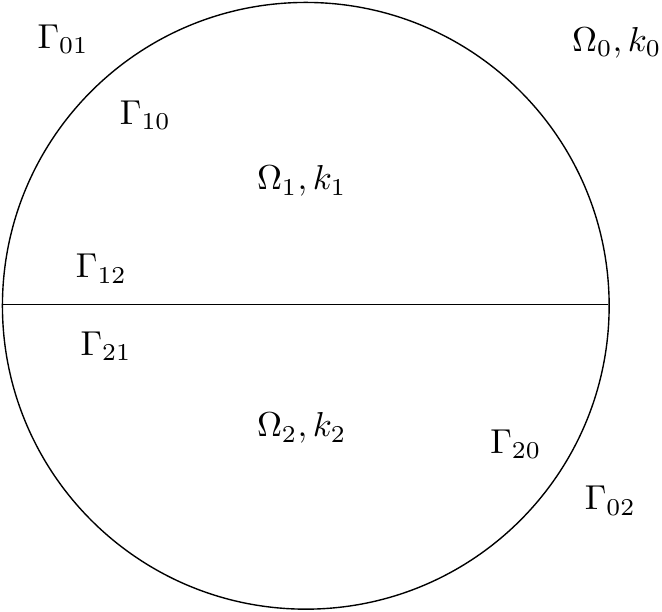}
\caption{Typical triple junction configuration.}
\label{fig:tj}
\end{figure}

We consider the problem of two dimensional electromagnetic scattering by structures that feature multiple junctions, i.e.~points where more than three interfaces of material discontinuity meet (e.g.~the structure displayed in Figure~\ref{fig:tj}). For the sake of presentation simplicity,  we focus our treatment of transmission problems with multiple junction domains on the two subdomain case depicted in Figure~\ref{fig:tj}. Specifically, we seek to solve the Helmholtz transmission problem that consists of finding the fields $u_0$, $u_1$, and $u_2$ solutions of the system of equations:
\begin{eqnarray}\label{system_tj}
  \Delta u_j +k_j^2 u_j &=&0\qquad {\rm in}\ \Omega_j,\\
  u_j+\delta_j^0 u^{inc}&=&u_\ell+\delta_\ell^0u^{inc}\qquad{\rm on}\ \Gamma_{j\ell}=\partial \Omega_j\cap \partial\Omega_\ell\nonumber,\\
  \varepsilon_j^{-1}(\partial_{n_j}u_j+\delta_j^0\partial_{n_j}u^{inc})&=&-\varepsilon_\ell^{-1}(\partial_{n_\ell}u_\ell+\delta_\ell^0\partial_{n_\ell}u^{inc})\qquad{\rm on}\ \Gamma_{j\ell}\nonumber,\\
 \lim_{r\to\infty}r^{1/2}(\partial u_0/\partial r-ik_0u_0)&=&0,\nonumber
\end{eqnarray}
where $\delta_j^0$ and $\delta_\ell^0$ stand for Kronecker operators, that is $\delta_j^0$ is the identity operator if $j=0$ and the null operator otherwise. Here, the wavenumber $k_j$ in the subdomain $\Omega_j$ is given by $k_j=\omega\sqrt{\varepsilon_j}$ in terms of the angular frequency $\omega>0$ and the electric permittivity $\varepsilon_j$. Throughout this paper  we assume that all the permittivities $\varepsilon_j$ are positive real numbers, extensions to more general cases being straightforward. The unit normal to the boundary $\partial\Omega_j$ is here denoted by $n_j$ and is assumed to point to the exterior of the subdomain $\Omega_j$. The incident field $u^{inc}$, on the other hand, is assumed to satisfy the Helmholtz equation with wavenumber $k_0$ in the unbounded domain $\Omega_0$. Finally, we point out that the well-posedness of the transmission problem~\eqref{system_tj}is discussed in~\cite{jerez-hanckes1} and references therein. 

In what follows, we review three main formulations of the transmission problem~\eqref{system_tj}. Two of them rely on boundary integral equations, and the third is a domain decomposition method. 

\section{Multitrace formulations\label{MS2}}

In this section we review the derivation of the MTF~\cite{jerez-hanckes1} of the transmission problem~\eqref{system_tj}. To this end, we make use of the four BIO associated with the Calder\'on calculus for a Lipschitz domain. Let~$D$ be a bounded domain in $\mathbb{R}^2$ whose boundary $\Gamma$ is a closed Lipschitz curve. Given a wavenumber~$k$, and a density $\varphi$ defined on $\Gamma$, we recall the definitions of the single layer potential~\cite{mclean:2000}
$$[SL_{\Gamma,k}(\varphi)](\mathbf{z}):=\int_\Gamma G_k(\mathbf{z}-\mathbf{y})\varphi(\mathbf{y})ds(\mathbf{y}),\ \mathbf{z}\in\mathbb{R}^2\setminus\Gamma,$$
and the double layer potential 
$$[DL_{\Gamma,k}(\varphi)](\mathbf{z}):=\int_\Gamma \frac{\partial G_k(\mathbf{z}-\mathbf{y})}{\partial\mathbf{n}(\mathbf{y})}\varphi(\mathbf{y})ds(\mathbf{y}),\ \mathbf{z}\in\mathbb{R}^2\setminus\Gamma,$$
where $G_k(\mathbf{x})=\frac{i}{4}H_0^{(1)}(k|\mathbf{x}|)$ represents the two-dimensional Green's function of the Helmholtz equation with wavenumber $k$, and $\mathbf{n}$ represents the unit normal pointing outside the domain $D$. Applying Dirichlet and Neumann exterior and interior traces on $\Gamma$--denoted by $\gamma_\Gamma^{D,1}$ and $\gamma_\Gamma^{D,2}$ and respectively $\gamma_\Gamma^{N,1}$ and $\gamma_\Gamma^{N,2}$--to the single and double layer potentials corresponding to the wavenumber $k$ and a density $\varphi$ we define the four Helmholtz BIO:
\begin{eqnarray}\label{traces}
\gamma_\Gamma^{D,1} SL_{\Gamma,k}(\varphi)&=&\gamma_\Gamma^{D,2} SL_{\Gamma,k}(\varphi)=S_{\Gamma,k}\varphi,\quad \gamma_\Gamma^{N,1} DL_{\Gamma,k}(\varphi)=\gamma_\Gamma^{N,2} DL_k(\varphi)=N_{\Gamma,k}\varphi\nonumber,\\
\gamma_\Gamma^{N,j} SL_{\Gamma,k}(\varphi)&=&(-1)^j\frac{\varphi}{2}+K_{\Gamma,k}^\top \varphi\quad j=1,2,\quad
\gamma_\Gamma^{D,j} DL_{\Gamma,k}(\varphi)=(-1)^{j}\frac{\varphi}{2}+K_{\Gamma,k}\varphi\quad j=1,2.\nonumber
\end{eqnarray}

\subsection{First kind (local) MTF}
The main idea in MTF is to use Calder\'on projectors in each subdomain $\Omega_j$, $j=0,1,2$ and to use continuity of Dirichlet and Neumann traces. Indeed, by applying Green's identities in a domain $\Omega_j$ we have that~\cite{KressColton}
\[
u_j=-DL_{\partial\Omega_j,j}u_{j} + SL_{\partial\Omega_j,j}\partial_{n_j}u_{j}{\quad{\mbox{in}\quad \Omega_j}}.
\]
For $j=0$, we first apply to the equation above the Dirichlet trace on $\partial\Omega_0$ to the domain $\Omega_0$ and obtain
\begin{equation}\label{eq:first_trace}
  \frac{1}{2}u_0 + K_{\partial\Omega_0,0}u_0 - \varepsilon_0S_{\partial\Omega_0,0}(\varepsilon_0^{-1}\partial_{n_0}u_0)=0\quad {\rm on}\quad \partial\Omega_0.
\end{equation}
The continuity conditions of Dirichlet traces on $\partial\Omega_0$, that is
\[
u_0=\begin{cases}
 u_1-u^{inc}, & {\rm on}\  \Gamma_{01}, \\
 u_2-u^{inc},  & {\rm on}\ \Gamma_{02},
\end{cases}
\]
can be expressed in the equivalent form
\begin{equation}\label{eq:cont_condD}
  u_0=X_{01}u_1 + X_{02}u_2 -u^{inc}, \quad{\rm on}\  \partial\Omega_0,
  \end{equation}
where the operators $X_{01}$ and $X_{02}$ are defined as
\[
X_{01} := E^{(0)}_{01}\Pi_{01} \quad X_{02} := E^{(0)}_{02}\Pi_{02}
\]
in terms of the restriction operators $\Pi_{ij}$ onto the interface $\Gamma_{ij}=\partial\Omega_i\cap\partial \Omega_j$--in the cases when $i\neq j$ and the domains $\Omega_i$ and $\Omega_j$ share an edge, and the extension-by-zero operators defined as
\[
E^{(0)}_{01}\psi_{01}:=\begin{cases}
 \psi_{01}, & {\rm on}\  \Gamma_{01} \\
 0,  & {\rm on}\ \Gamma_{02}
\end{cases}
\qquad
E^{(0)}_{02}\psi_{02}:=\begin{cases}
 0, & {\rm on}\  \Gamma_{01} \\
 \psi_{02},  & {\rm on}\ \Gamma_{02}
\end{cases}
\]
for two functions $\psi_{01}$ and $\psi_{02}$ defined on $\Gamma_{01}$ and $\Gamma_{02}$ respectively. Using the reformulation~\eqref{eq:cont_condD} we can re-express equation~\eqref{eq:first_trace} in the following form
\begin{equation}\label{eq:first_trace_X}
 K_{\partial\Omega_0,0}u_0 - \varepsilon_0S_{\partial\Omega_0,0}(\varepsilon_0^{-1}\partial_{n_0}u_0)+\frac{1}{2}X_{01}u_1+\frac{1}{2}X_{02}u_2=\frac{1}{2}u^{inc}\quad {\rm on}\ \partial\Omega_0.
\end{equation}
\begin{remark}
  The expressions $u_1$ and $u_2$ in equation~\eqref{eq:first_trace_X} should be understood as Dirichlet traces $u_1|_{\partial\Omega_1}$ and $u_2|_{\partial\Omega_2}$. We prefer the simpler notation as the meaning of those expressions should be grasped from the domain of definition of the operators that act upon them. Furthermore, all the restriction/extension operators can be defined for distributions~\cite{jerez-hanckes1}.
\end{remark}
Applying the normal derivative with respect to $\mathbf{n}_0$ to the same Green's identities above we obtain
\begin{equation}\label{eq:second_trace}
  \frac{1}{2}\varepsilon_0^{-1}\partial_{n_0}u_0 - K_{\partial\Omega_0,0}^\top(\varepsilon_0^{-1}\partial_{n_0}u_0) + \varepsilon_0^{-1}N_{\partial\Omega_0,0}u_0=0\quad {\rm on}\ \partial\Omega_0.
\end{equation}
The continuity conditions of Neumann traces on $\partial\Omega_0$ 
\[
\varepsilon_0^{-1}\partial_{n_0}u_0:=\begin{cases}
 -\varepsilon_1^{-1}\partial_{n_1}u_1-\varepsilon_0^{-1}\partial_{n_0}u^{inc}, & {\rm on}\  \Gamma_{01}, \\
 -\varepsilon_2^{-1}\partial_{n_2}u_2-\varepsilon_0^{-1}\partial_{n_0}u^{inc},  & {\rm on}\ \Gamma_{02},
\end{cases}
\]
can be expressed in the equivalent form
\begin{equation}\label{eq:cont_condN}
  \varepsilon_0^{-1}\partial_{n_0}u_0=-X_{01}(\varepsilon_1^{-1}\partial_{n_1}u_1) - X_{02}(\varepsilon_2^{-1}\partial_{n_2}u_2) -\varepsilon_0^{-1}\partial_{n_0}u^{inc}, \quad{\rm on}\  \partial\Omega_0,
  \end{equation}
and thus, equation~\eqref{eq:second_trace} can be expressed as 
\begin{equation}\label{eq:second_traceX}
  \varepsilon_0^{-1}N_{\partial\Omega_0,0}u_0 - K_{\partial\Omega_0,0}^\top(\varepsilon_0^{-1}\partial_{n_0}u_0) -\frac{1}{2}X_{01}(\varepsilon_1^{-1}\partial_{n_1}u_1) - \frac{1}{2}X_{02}(\varepsilon_2^{-1}\partial_{n_2}u_2)=\frac{1}{2}\varepsilon_0^{-1}\partial_{n_0}u^{inc}\quad {\rm on}\ \partial\Omega_0.
\end{equation}
We define the Calder\'on block operators $\mathcal{C}_j$ satisfying 
\[
\mathcal{C}_j:=\begin{bmatrix}K_{\partial\Omega_j,j} & -\varepsilon_jS_{\partial\Omega_j,j}\\
\varepsilon_j^{-1}N_{\partial\Omega_j,j} & - K_{\partial\Omega_j,j}^\top
\end{bmatrix},\ j=0,1,2,\qquad \mathcal{C}_j^2=\frac{1}{4}\begin{bmatrix}I & 0\\0 & I\end{bmatrix},
\]
and also introduce the following matrix operators:
\[
\mathbb{X}_{01}=\begin{bmatrix}X_{01} & 0\\ 0 & -X_{01}\end{bmatrix}\qquad \mathbb{X}_{02}=\begin{bmatrix}X_{02} & 0\\ 0 & -X_{02}\end{bmatrix}.
\]
Denoting the Cauchy data on $\partial\Omega_j$ by
\[
\mathbf{u}_j:=\begin{bmatrix}u_j\\ \varepsilon_j^{-1}\partial_{n_j}u_j\end{bmatrix},\ j=0,1,2,
\]
and the incoming data $\mathbf{u}^{inc}$  on $\partial\Omega_0$ by
\[
\mathbf{u}^{inc}:=\begin{bmatrix}u^{inc}\\ \varepsilon_0^{-1}\partial_{n_0}u^{inc}\end{bmatrix},
\]
we can re-write equations~\eqref{eq:first_trace_X} and~\eqref{eq:second_traceX} in the form
\begin{equation}\label{f_eq_m}
  \mathcal{C}_0\mathbf{u}_0+\frac{1}{2}\mathbb{X}_{01}\mathbf{u}_1+\frac{1}{2}\mathbb{X}_{02}\mathbf{u}_2=\frac{1}{2}\mathbf{u}^{inc}\quad{\rm on}\ \partial\Omega_0.
\end{equation}
Similar techniques applied to $u_1$ in the domain $\Omega_1$ and respectively $u_2$ in the domain $\Omega_2$ lead to two additional equations in the vein of equation~\eqref{f_eq_m} leading to the local MTF~\cite{jerez-hanckes1}:
\begin{equation}\label{eq:MTF1}
  \begin{bmatrix}
    \mathcal{C}_0& \mathbb{X}_{01} & \mathbb{X}_{02}\\
    \mathbb{X}_{10} & \mathcal{C}_1 & \mathbb{X}_{12}\\
    \mathbb{X}_{20} & \mathbb{X}_{21} & \mathcal{C}_2
  \end{bmatrix}
    \begin{bmatrix}\mathbf{u}_0\\ \mathbf{u}_1\\ \mathbf{u}_2\end{bmatrix}=\frac{1}{2}\begin{bmatrix}\mathbf{u}^{inc}\\ -\mathbb{X}_{10}\mathbf{u}^{inc} \\ -\mathbb{X}_{20}\mathbf{u}^{inc}\end{bmatrix},
\end{equation}
where the additional matrix operators $\mathbb{X}_{j\ell}$ are similarly defined. The well-posedness of the local MTF in the multitrace space $\prod_{j=0}^2 \left(H^{1/2}(\partial \Omega_j)\times  H^{-1/2}(\partial \Omega_j)\right)$ was established in the same reference~\cite{jerez-hanckes1}. We note that  owing to the presence of both weakly singular and hypersingular operators in the Calder\'on projectors $\mathcal{C}_j, j=0,1,2$, MTF~\eqref{eq:MTF1} are boundary integral equations of the first kind. As such, numerical solution of the MTF typically resort to preconditioning, an issue we return to in Section~\ref{num}.

\subsection{Second kind (global) MTF}
\label{sk}

An alternative boundary integral equation formulation of the transmission problem~\eqref{system_tj} can be derived using layer potentials defined on the skeleton $\Gamma:=\Gamma_{01}\cup\Gamma_{12}\cup\Gamma_{02}$ whose segments are oriented clockwise~\cite{jerez-hanckes2}. In this approach, the solutions $u_j,j=0,1,2$ of the transmission problem~\eqref{system_tj} are sought in the form
\begin{equation}\label{eq:layer_tj}
  u_j(\mathbf{x}):=SL_{\Gamma,j}\ v-\varepsilon_j DL_{\Gamma,j}\ p,\quad \mathbf{x}\in\Omega_j
\end{equation}
where $v$ and $p$ are densities defined on the skeleton $\Gamma$ and the double layer operators are defined with respect to exterior unit normals $\mathbf{n}$ with respect to the orientation chosen on $\Gamma$. Here we used the same convention that the index $j$ in equation~\eqref{eq:layer_tj} refers to the wavenumber $k_j$ for $j=0,1,2$. The enforcement of the continuity conditions of Dirichlet and scaled Neumann traces on the interfaces between subdomains leads to the following system of boundary integral equations, for all $\ell,j=0,1,2$ such that $\ell<j$:
\begin{eqnarray}\label{eq:sk_system}
  \frac{\varepsilon_\ell+\varepsilon_j}{2}p+(-1)^\ell[S_\ell-S_j]v+(-1)^\ell[\varepsilon_\ell K_\ell-\varepsilon_j K_j]p&=&\begin{cases}
 -u^{inc},\ {\rm if}\ \ell=0\\
 0\ {\rm otherwise}
  \end{cases}\quad {\rm on}\ \Gamma_{\ell j},\ \nonumber\\
  \frac{\varepsilon_\ell^{-1}+\varepsilon_j^{-1}}{2}v+(-1)^\ell[\varepsilon_j^{-1}S_j-\varepsilon_\ell^{-1}S_\ell]v+(-1)^\ell[N_j-N_\ell]p&=&\begin{cases}
 -\varepsilon_0^{-1}\partial_nu^{inc},\ {\rm if}\ \ell=0\\
 0\ {\rm otherwise}
  \end{cases}\ {\rm on}\ \Gamma_{\ell j}.\nonumber\\
\end{eqnarray}
In what follows we refer to the formulation~\eqref{eq:sk_system} by the acronym g-MTF. All of the boundary integral operators featured in equation~\eqref{eq:sk_system} are defined on the skeleton $\Gamma$. The presence of (a) multiples of the densities $v$ and $p$ as well as (b) differences of hypersingular operators $N_j-N_{\ell}$ suggest referring to equations~\eqref{eq:sk_system} as second kind boundary integral equations, hence the acronym CFIESK. We note that in contast to the MTF~\eqref{eq:MTF1}, the second kind equations~\eqref{eq:sk_system} require one (scalar) unknown per interface $\Gamma_{j\ell}$. To the best of our knowledge, the well posedness of second kind boundary integral equations~\eqref{eq:sk_system} is still an open question, see the discussion on this in~\cite{claeys2015second}. Our numerical experiments indicate that these formulations appear to be well posed for a wide range of wavenumbers $k_j,j=0,1,2$, thus corroborating the findings in~\cite{claeys2015second}.

\section{Domain decomposition approach with classical Robin boundary conditions\label{MS}}

DDM are natural candidates for numerical solution of transmission problems~\eqref{system_tj}. A non-overlapping domain decomposition approach for the solution of equations~\eqref{system_tj} consists of solving subdomain problems with matching Robin boundary conditions on the common subdomain interfaces~\cite{Depres}. Indeed, this procedure amounts to computing the subdomain solutions:
\begin{eqnarray}\label{DDM}
  \Delta u_j +k_j^2 u_j &=&0\qquad {\rm in}\ \Omega_j,\\
  \varepsilon_j^{-1}(\partial_{n_j}u_j+\delta_j^0\partial_{n_j}u^{inc})+i\eta (u_j+\delta_j^0 u^{inc})&=&-\varepsilon_\ell^{-1}(\partial_{n_\ell}u_\ell+\delta_\ell^0 \partial_{n_\ell}u^{inc})+i\eta (u_\ell+\delta_\ell^0 u^{inc})\qquad{\rm on}\ \Gamma_{j\ell}\nonumber.
\end{eqnarray}
In equations~\eqref{DDM}, $\eta$ is assumed to be a positive number. The latter requirement is needed to ensure the well-posedness of the impedance boundary value Helmholtz problem in the exterior domain $\Omega_0$~\cite{KressColton}.\footnote{In all the numerical examples in Section~\ref{num} we took $\eta=k_0$.}

In order to describe the DDM method more concisely we introduce subdomain Robin-to-Robin (RtR) maps~\cite{Collino1}. Given a subdomain $\Omega_j$ we define the RtR map $\mathcal{S}^j$ in the following manner:
\begin{equation}\label{RtRboxj}
   \mathcal{S}^j(\psi_j):=(\varepsilon_j^{-1}\partial_{n_j} u_j-i\eta\ u_j)|_{\partial \Omega_j},\quad j=0,1,2,
 \end{equation}
 where $u_j$ is the solution of the following problem:
 \begin{eqnarray*}
   \Delta u_j+k_j^2u_j&=&0\ {\rm in}\ \Omega_j,\\
   \varepsilon_j^{-1}\partial_{n_j} u_j+i\eta u_j&=&\psi_j\ {\rm on}\ \partial \Omega_j.
 \end{eqnarray*}
In the case when $\Omega_j$ is the exterior domain $\Omega_0$, we further require in the definition of the RtR map $\mathcal{S}^0$ that $u_0$ is radiative at infinity. The DDM method computes the global Robin data $f=[f_1\ f_2\ f_0]^\top$ --the reader should observe the coordinate order-- with
\[
f_{j}:=(\varepsilon_j^{-1}\partial_{n_j}u_j+i\eta\ u_j)|_{\partial \Omega_j},\ j=0,1,2
\]
as the solution of the following linear system that incorporates the subdomain RtR maps $\mathcal{S}^j,j=0,1,2$, previously defined
 \begin{equation}\label{ddm}
 (I+A)f=g,\quad A:=X\mathcal{S},\quad \mathcal{S}=\begin{bmatrix}\mathcal{S}^1&0&0\\0&\mathcal{S}^2&0\\0&0&\mathcal{S}^0\end{bmatrix}, \quad  X:=\begin{bmatrix}0 & X_{12}& X_{10}\\ X_{21}&0& X_{20}\\ X_{01}& X_{02}&0\end{bmatrix}.
 \end{equation}
with right-hand side $g=[g_1\ g_2\ g_0]^\top$ wherein
 \begin{eqnarray}\label{rhs_ddm}
   g_1&=& X_{01}(-\varepsilon_0^{-1}\partial_{n_0}u^{inc}+i\eta u^{inc})|_{\partial\Omega_0}\nonumber\\
   g_2&=&X_{02}(-\varepsilon_0^{-1}\partial_{n_0}u^{inc}+i\eta u^{inc})|_{\partial\Omega_0}\nonumber\\
   g_0&=&(-\varepsilon_0^{-1}\partial_{n_0}u^{inc}-i\eta u^{inc})|_{\partial\Omega_0}.\nonumber
   \end{eqnarray}
 \begin{remark}
   The domains $\Omega_j$, $j\neq 0$, can be further subdivided into smaller subdomains, in which case the DDM system~\eqref{ddm} has to be augmented to incorporate the additional Robin data on the new interfaces. The size of the subdomains--in terms of wavelengths--should ideally be such that the computation/application of the corresponding RtR operators can be performed efficiently.  
   \end{remark}
 We note that the matrix $A$ in equation~\eqref{ddm} is not stored in practice, and, due to its possibly large size, the DDM linear system~\eqref{ddm} is typically solved in practice via iterative methods. Iterative solvers (e.g.~Jacobi, GMRES) for the solution of DDM linear systems of the type described in equation~\eqref{ddm} require large numbers of iterations, especially in the case of larger numbers of subdomains (e.g.~see Table~\ref{comp5} in Section~\ref{num}). This shortcoming can be attributed to the choice of classical Robin boundary conditions and the outflow/inflow of information from a subdomain to its neighboring subdomains associated with it. Ideally, the subdomain boundary conditions have to be chosen so that information flows out of the subdomain and no information is reflected back into the subdomain. This can be achieved if the term $i\eta$ in equations~\eqref{ddm} is replaced by the adjacent subdomain Dirichlet-to-Neumann (DtN) operator restricted to the common interface---in this way the Jacobi scheme converges in precisely two iterations~\cite{Nataf,HJP13}. Since DtN maps are not always well defined and expensive to calculate even when properly defined, easily computable approximations of DtN maps can be employed effectively to lead to faster convergence rates of GMRES solvers for DDM algorithms~\cite{boubendirDDM}. We present in Section~\ref{GIBC} a non-overlapping DDM based on square root approximations of DtN operators.

\begin{figure}
\centering
\includegraphics[height=60mm]{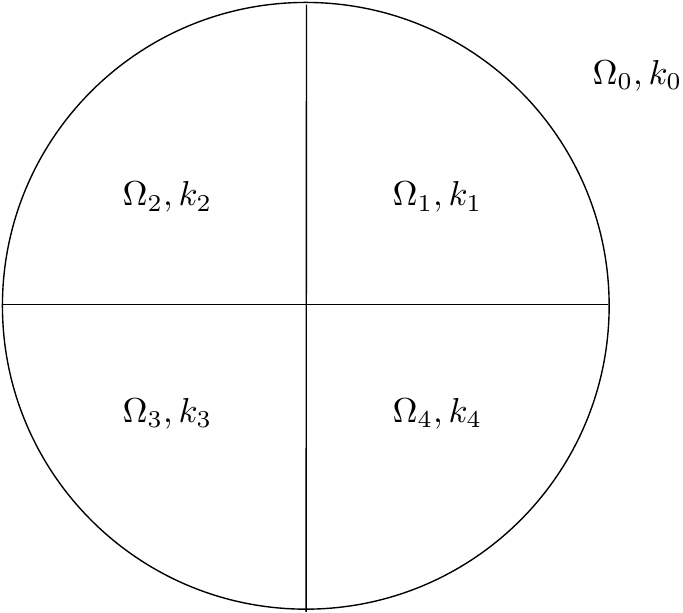}
\caption{Five subdomain configuration.}
\label{fig:4sbds}
\end{figure}

In order to bypass the shortcomings of iterative solvers for DDM formulations~\eqref{ddm}, we resort to a direct solver that consists of application of recursive Schur complements. We begin by further splitting the subdomain data $f$ into interface components 
\[
f_{j\ell}:=(\varepsilon_j^{-1}\partial_{n_j}u_j+i\eta\ u_j)|_{\partial\Omega_j\cap\partial \Omega_\ell},\ 0\leq j,\ell\leq 2,\ \operatorname{meas}(\partial \Omega_\ell\cap\partial \Omega_j)\neq 0,
\]
 and by rewriting the linear sytem~\eqref{ddm} in terms of the new unknown arrangement. To this end, we further decompose the RtR operators $\mathcal{S}^1$ and $\mathcal{S}^2$ in the block form
\[
\mathcal{S}^{1}=\begin{bmatrix}\mathcal{S}^1_{01,10}&\mathcal{S}^1_{01,12}\\ \mathcal{S}^1_{21,10}& \mathcal{S}^1_{21,12}\end{bmatrix}\qquad\text{and}\qquad \mathcal{S}^{2}=\begin{bmatrix}\mathcal{S}^{2}_{02,20}&\mathcal{S}^{2}_{02,21}\\ \mathcal{S}^{2}_{12,20}& \mathcal{S}^{2}_{12,21}\end{bmatrix},
\]
where the block operators are defined informally as
\begin{equation}\label{SL}
  \begin{bmatrix}\mathcal{S}^1_{01,10}&\mathcal{S}^1_{01,12}\\ \mathcal{S}^1_{21,10}& \mathcal{S}^1_{21,12}\end{bmatrix}\begin{bmatrix}(\varepsilon_1^{-1}\partial_{n_1} u_1 +i\eta u_1)|_{\Gamma_{10}}\\(\varepsilon_1^{-1}\partial_{n_1} u_1 +i\eta u_1)|_{\Gamma_{12}}\end{bmatrix}=\begin{bmatrix}(\varepsilon_1^{-1}\partial_{n_1} u_1 - i\eta u_1)|_{\Gamma_{10}}\\(\varepsilon_1^{-1}\partial_{n_1} u_1 - i\eta u_1)|_{\Gamma_{12}}\end{bmatrix}
\end{equation}
and
\begin{equation}\label{SR}
\begin{bmatrix}\mathcal{S}^{2}_{02,20}&\mathcal{S}^{2}_{02,21}\\ \mathcal{S}^{2}_{12,20}& \mathcal{S}^{2}_{12,21}\end{bmatrix}\begin{bmatrix}(\varepsilon_2^{-1}\partial_{n_2} u_{2} +i\eta u_{2})|_{\Gamma_{20}}\\(\varepsilon_2^{-1}\partial_{n_2} u_{2} +i\eta u_{2})|_{\Gamma_{21}}\end{bmatrix}=\begin{bmatrix}(\varepsilon_2^{-1}\partial_{n_2} u_{2} - i\eta u_{2})|_{\Gamma_{20}}\\(\varepsilon_2^{-1}\partial_{n_2} u_{2} - i\eta u_{2})|_{\Gamma_{21}}\end{bmatrix}.
\end{equation}
We similalry split the operator $\mathcal{S}^0=[\mathcal{S}^0_{10}\quad \mathcal{S}^0_{20}]$. These block decompositions/splittings of the RtR maps $\mathcal{S}^j,j=0,1,2$ allow us to write explicitly the DDM system~\eqref{ddm} in the form:
\begin{equation}\label{DDM_explicit}
  \begin{bmatrix}
    I & \mathcal{S}^2_{12,21} & 0 & \mathcal{S}^2_{12,20} & 0\\
    \mathcal{S}^1_{21,12} & I & \mathcal{S}^1_{21,10} & 0 & 0\\
    0 & 0 & I & 0 & \mathcal{S}^0_{10}\\
    0 & 0 & 0 & I & \mathcal{S}^0_{20}\\
    X_{01}\mathcal{S}^1_{01,12} & X_{02}\mathcal{S}^2_{02,21} & X_{01}\mathcal{S}^1_{01,10} & X_{02}\mathcal{S}^2_{02,20} & I\\
  \end{bmatrix}
  \begin{bmatrix}
    f_{12}\\f_{21}\\f_{10}\\f_{20}\\f_0
  \end{bmatrix}=
  \begin{bmatrix}
    0\\0\\(-\varepsilon_0^{-1}\partial_{n_0}u^{inc}+i\eta u^{inc})|_{\Gamma_{01}}\\(-\varepsilon_0^{-1}\partial_{n_0}u^{inc}+i\eta u^{inc})|_{\Gamma_{02}}\\(-\varepsilon_0^{-1}\partial_{n_0}u^{inc}-i\eta u^{inc})|_{\Gamma_0}
  \end{bmatrix}.
\end{equation}
We note that the matrix in the new linear system~\eqref{DDM_explicit} exhibits a block sparsity pattern that resembles that of matrices corresponding to two-dimensional finite difference operators that  use five-point Laplaceans. We use a Schur complement approach for the solution of the linear system~\eqref{DDM_explicit} that consists of elimination of the Robin data $(f_{12},f_{21})$ corresponding to the interior interface $\Gamma_{12}$ from the linear system~\eqref{DDM_explicit}. This procedure requires using the inverse of the matrix
\[
\mathcal{D}_{12}:=\begin{bmatrix} I & \mathcal{S}^2_{12,21} \\ \mathcal{S}^1_{21,12} & I \end{bmatrix}
\]
which can be computed explicitly
\begin{equation}\label{inv_matrix_explicit}
  \mathcal{D}_{12}^{-1}=\begin{bmatrix}
I+\mathcal{S}^2_{12,21}(I-\mathcal{S}^1_{21,12}\mathcal{S}^2_{12,21})^{-1}\mathcal{S}^1_{21,12} & -\mathcal{S}^2_{12,21}(I-\mathcal{S}^1_{21,12}\mathcal{S}^2_{12,21})^{-1}\\
-(I-\mathcal{S}^1_{21 ,12}\mathcal{S}^2_{12,21})^{-1}\mathcal{S}^1_{21,12} &
(I-\mathcal{S}^1_{21,12}\mathcal{S}^2_{12,21})^{-1}
\end{bmatrix}.
\end{equation}
Using the Schur complement in equation~\eqref{DDM_explicit} we obtain the following linear system
\begin{equation}\label{eq:reduced}
  \begin{bmatrix}I & \mathcal{S}^0\\\mathcal{S}^{int}& I \end{bmatrix}\begin{bmatrix}f_0^{int}\\f_0\end{bmatrix}=\begin{bmatrix}(-\varepsilon_0^{-1}\partial_{n_0}u^{inc}+i\eta u^{inc})|_{\Gamma_{0}}\\(-\varepsilon_0^{-1}\partial_{n_0}u^{inc}-i\eta u^{inc})|_{\Gamma_{0}}\end{bmatrix},
  \end{equation}
where $f_0^{int}=[f_{10}\quad f_{20}]^\top$. Interestingly, the operator $\mathcal{S}^{int}$ that features in equation~\eqref{eq:reduced} is itself an RtR operator defined as 
\begin{equation}\label{RtRint}
   \mathcal{S}^{int}(\psi):=(-\partial_{n_0} u-i\eta\ \varepsilon(x)\ u)|_{\partial \Omega_0},\quad \varepsilon(x)=\varepsilon_1,\ x\in\Gamma_{10},\ \varepsilon(x)=\varepsilon_2,\ x\in \Gamma_{20},
 \end{equation}
 where $u$ is the solution of the following problem:
 \begin{eqnarray*}
   \Delta u+k(x)^2u&=&0\ {\rm in}\ \Omega_1\cup\Omega_2,\quad k(x)=k_1,\ x\in\Omega_1,\ k(x)=k_2,\ x\in\Omega_2,\\
   -\partial_{n_0} u+i\eta\ \varepsilon(x) u&=&\psi\ {\rm on}\ \partial \Omega_0,
 \end{eqnarray*}
 and $u$ and $\varepsilon(x)^{-1}\partial_{n}u$ are continuous across $\Gamma_{12}$ where $\varepsilon(x)=\varepsilon_1,\ x\in\Omega_1,\ \varepsilon(x)=\varepsilon_2,\ x\in\Omega_2$, and $n=n_1=-n_2$ on $\Gamma_{12}$. For obvious reasons, the procedure by which $\mathcal{S}^{int}$ was obtained above can be interpreted as the merging of subdomain RtR maps $\mathcal{S}^j,j=1,2$.
 \begin{remark}
   In the case of five subdomain configuration depicted in Figure~\ref{fig:4sbds}, the Schur complement elimination procedure proceeds by rearranging the subdomain Robin data in the order
   \[
   f=[f_{12}\ f_{21}\ f_{34}\ f_{43}\ f_{23}\ f_{14}\ f_{32}\ f_{41}\ f_{10}\ f_{20}\ f_{30}\ f_{40}\ f_0]^\top
   \]
   and then (1) first eliminating the unknown pairs $[f_{12},\ f_{21}]$ and $[f_{34},\ f_{43}]$ in the first stage, and (2) then subsequently eliminating the unknown pair $f_{top}:=[f_{23},\ f_{14}]$ and $f_{bottom}:=[f_{32},\ f_{41}]$. After performing the elimination stage we arrive at a sytem of the form~\eqref{eq:reduced} with $f_0^{int}=[f_{10}\ f_{20}\ f_{30}\ f_{40}]$ and the same interpretation of the operator $\mathcal{S}^{int}$ in the domain $\Omega_1\cup\Omega_2\cup\Omega_3\cup\Omega_4$. 
   \end{remark}
The system~\eqref{eq:reduced} reduces in turn to the following equation involving the Robin unknown $f_0$ defined on the other interface $\partial\Omega_0$:
 \begin{equation}\label{eq:final}
   (I-\mathcal{S}^{int}\mathcal{S}^0)f_0=(-\varepsilon_0^{-1}\partial_{n_0}u^{inc}-i\eta u^{inc})|_{\Gamma_{0}}-\mathcal{S}^{int}(-\varepsilon_0^{-1}\partial_{n_0}u^{inc}+i\eta u^{inc})|_{\Gamma_{0}}.
 \end{equation}
 Ideally, the linear system~\eqref{eq:final} should be solved by direct methods as well. Alternatively, if size of the exterior Robin unknown $f_0$ becomes too large for application of direct linear algebra solvers, Krylov subspace iterative solvers such as GMRES~\cite{SaadSchultz} can be employed for the solution of the reduced system~\eqref{eq:final}. We note that the operator $I-\mathcal{S}^{int}\mathcal{S}^0$ is a pseudodifferential operator of order $-1$ (this can be rigorously established in the case when $k_1=k_2$), and thus equation~\eqref{eq:final} does not possess ideal spectral properties for fast convergence of iterative Krylov subspace solvers. Alternative interior/exterior coupling strategies that possess superior spectral properties is subject of ongoing investigations.   

 \subsection{Analysis of the Schur complement DDM algorithm\label{inv}}

 As it can be seen from the description in Section~\ref{MS}, the DDM Schur complement procedure hinges on the invertibility of the operators $I-\mathcal{S}^1_{21,12}\mathcal{S}^2_{12,21}$, an issue which we examine in what follows. A more precise definition of the block operators in equations~\eqref{SL} and~\eqref{SR} can be given by considering the partial Robin-to-Robin maps:
\[
\mathcal{S}^j_\ell(0,\varphi_{j\ell}) := (\varepsilon_j^{-1}\partial_{n_j} w_j+ i\eta w_j)|_{\partial \Omega_j},\ j,\ell\geq1,\ j\neq \ell
\]
where
\begin{eqnarray}\label{Robin_int}
  \Delta w_j+k_j^2 w_j&=&0\ {\rm in}\  \Omega_j, \nonumber\\
  \varepsilon_j^{-1}\partial_{n_j} w_j-i\eta w_j&=&\varphi_{j\ell}\ {\rm on}\quad \Gamma_{j\ell},\\
  \varepsilon_j^{-1}\partial_{n_j} w_j-i\eta w_j&=&0\ {\rm on}\quad \Gamma_{j0}.\nonumber
\end{eqnarray}
Then we define
\[
\mathcal{S}^{1}_{21,10}\varphi_{12}\:=  (\varepsilon_1^{-1}\partial_{n_1} w_1+i\eta w_1)|_{\Gamma_{10}},\quad \mathcal{S}^{1}_{21,12}\varphi_{12}\:=  (\varepsilon_1^{-1}\partial_{n_1} w_1+i\eta w_1)|_{\Gamma_{12}},
\]
and similarly for the other block operators.

The next results sheds light into the mapping properties of the operators $\mathcal{S}^1_{21,12}$ and $\mathcal{S}^2_{12,21}$:
\begin{theorem}\label{spectral_radius}
  The following inequalities hold true
  \[
  \|\mathcal{S}^1_{21,12}\|_{L^2(\Gamma_{12})\to L^2(\Gamma_{12})}\leq 1\quad \|\mathcal{S}^2_{12,21}\|_{L^2(\Gamma_{12})\to L^2(\Gamma_{12})}\leq 1.
  \]
\end{theorem}
\begin{proof} This follow easily from the definition of the block operators $\mathcal{S}^1_{21,12}$ and $\mathcal{S}^2_{12,21}$ and the fact that $\mathcal{S}^1_2$ and $\mathcal{S}^2_1$ are isometries in $L^2(\partial \Omega_1)$ and $L^2(\partial\Omega_2)$ respectively.
\end{proof}

We establish next the injectivity of the operator $I-\mathcal{S}^1_{21,12}\mathcal{S}^2_{12,21}$:
\begin{theorem}\label{thm_inv}
  The operator $I-\mathcal{S}^1_{21,12}\mathcal{S}^2_{12,21}:L^2(\Gamma_{12})\to L^2(\Gamma_{12})$ is injective.
\end{theorem}
\begin{proof}  Let $\psi\in \operatorname{Ker}(I-\mathcal{S}^1_{21,12}\mathcal{S}^2_{12,21})$, that is
 \[
  \mathcal{S}^1_{21,12}\mathcal{S}^2_{12,21}\psi=\psi.
  \]
  It is immediate to see that the operator $\mathcal{S}^1_{21,12}$ is invertible, and its inverse is defined as
  \[
  (\mathcal{S}^1_{21,12})^{-1}=\Pi_{12} \widetilde{\mathcal{S}}^{1}_{2}
  \]
  where
  \[
\widetilde{\mathcal{S}}^{1}_2(0,\psi_{12}) := (\varepsilon_1^{-1}\partial_{n_1} u- i\eta u)|_{\partial \Omega_1},
\]
with $u$ being the solution of the following boundary value problem:
\begin{eqnarray*}
  \Delta u+k_1^2 u&=&0\ {\rm in}\  \Omega_1\\
  \varepsilon_1^{-1}\partial_{n_1} u+i\eta u&=&\psi_{12}\ {\rm on}\quad \Gamma_{12},\\
  \varepsilon_1^{-1}\partial_{n_1} u-i\eta u&=&0\ {\rm on}\quad \Gamma_{10}.
\end{eqnarray*}
Define then $v_2$ to be the unique solution of 
\begin{eqnarray*}
  \Delta v_2+k_2^2 v_2&=&0\ {\rm in}\  \Omega_2\\
  \varepsilon_2^{-1}\partial_{n_2} v_2-i\eta v_2&=&\psi\ {\rm on}\quad \Gamma_{21},\\
  \varepsilon_2^{-1}\partial_{n_2} v_2-i\eta v_2&=&0\ {\rm on}\quad \Gamma_{20},
\end{eqnarray*}
and $v_1$ to be the unique solution of
\begin{eqnarray*}
  \Delta v_1+k_1^2 v_1&=&0\ {\rm in}\  \Omega_1\\
  \varepsilon_1^{-1}\partial_{n_1} v_1+i\eta v_1&=&\psi\ {\rm on}\quad \Gamma_{12},\\
  \varepsilon_1^{-1}\partial_{n_1} v_1-i\eta v_1&=&0\ {\rm on}\quad \Gamma_{10}.
\end{eqnarray*}
The relation $\mathcal{S}^1_{21,12}\mathcal{S}^2_{12,21}\psi=\psi$, which is equivalent to
$\mathcal{S}^2_{12,21}\psi=(\mathcal{S}^1_{21,12})^{-1}\psi$ implies that
\[
\varepsilon_1^{-1}\partial_{n_1} v_1-i\eta v_1=\varepsilon_2^{-1}\partial_{n_2} v_2+i\eta v_2\quad {\rm on}\quad \Gamma_{12}.
\]
The last equation together with the fact that
\[
\varepsilon_1^{-1}\partial_{n_1} v_1+i\eta v_1=\varepsilon_2^{-1}\partial_{n_2} v_2-i\eta v_2=\psi\quad {\rm on}\quad \Gamma_{12}
\]
imply that
\begin{equation}\label{equality}
  v_1 = -v_2, \quad \varepsilon_1^{-1}\partial_{n_1} v_1=\varepsilon_2^{-1}\partial_{n_2} v_2\quad {\rm on}\quad \Gamma_{12}.
  \end{equation}
Defining then
\[
v:=\left\{ \begin{array}{cc} v_1 &\quad{\rm in} \quad \Omega_1,\\ -v_2 &\quad {\rm in }\quad \Omega_2,\end{array}\right.
  \]
  we see from equation~\eqref{equality} that $v$ is continuous across $\Gamma_{12}$, and, if we take into account that $n:=n_1=-n_{2}$ on $\Gamma_{12}$, so is $\varepsilon_j^{-1}\partial_{n}v$ across $\Gamma_{12}$. In addition, $v$ satisfies the following boundary value problem
  \begin{eqnarray*}
  \Delta v+k(x)^2 v&=&0\ {\rm in}\  \Omega_1\cup\Omega_2,\quad k(x)=k_1,\ x\in\Omega_1,\ k(x)=k_2,\ x\in\Omega_2\\
  \varepsilon(x)^{-1}\partial_{n} v-i\eta v&=&0\ {\rm on}\quad \partial(\Omega_1\cup \Omega_2),\quad \varepsilon(x)=\varepsilon_1,\ x\in\Omega_1,\ \varepsilon(x)=\varepsilon_2,\ x\in\Omega_2.
  \end{eqnarray*}
Taking into account the continuity properties of $v$ and its normal derivative across $\Gamma_{12}$, it follows immediately that
  \[
  \int_{\Omega_1\cup\Omega_2}(-\omega^2 |v|^2 + \varepsilon(x)^{-1}|\nabla v|^2)dx=-i\eta\ \int_{\partial(\Omega_1\cup\Omega_2)}|v|^2\ ds,
  \]
  which implies that $v=0$ on $\partial(\Omega_1\cup\Omega_2)$. Given the zero impedance boundary condition on $\partial(\Omega_1\cup\Omega_2)$, we also have that $\partial_nv=0$ on $\partial(\Omega_1\cup\Omega_2)$, and hence $v=0$ in $\Omega_1\cup\Omega_2$. The last relation, in turn,  implies that $\psi=0$ on $\Gamma_{12}$ which completes the proof.
\end{proof}

The results in Theorem~\ref{spectral_radius} and Theorem~\ref{thm_inv} do not guarantee the invertibility of the operator $I-\mathcal{S}^1_{21,12}\mathcal{S}^2_{12,21}$ in $L^2(\Gamma_{12})$. Indeed, the operators $\mathcal{S}^1_{21,12}$ and $\mathcal{S}^2_{12,21}$ are not strict contractions, nor are they compact in the aforementioned spaces. Nevertheless, one could attempt to replace the inverses of the operators $I-\mathcal{S}^1_{21,12}\mathcal{S}^2_{12,21}$--and its counterparts--needed in the DDM Schur complement algorithm by their truncated Neumann series expansions:
\begin{equation}\label{eq:Neumannseries}
  (I-\mathcal{S}^1_{21,12}\mathcal{S}^2_{12,21})^{-1}\approx \sum_{m=0}^{N_{trunc}} (\mathcal{S}^1_{21,12}\mathcal{S}^2_{12,21})^m.
\end{equation}
Given that for small values of the parameter $N_{trunc}$  the evaluation of the expression in the right-hand side of equation~\eqref{eq:Neumannseries} is less expensive than the full inversion, we explore in Section~\ref{num} the effectiveness of this approach.

\subsection{DDM with generalized Robin boundary conditions~\label{GIBC}}

The rate of convergence of iterative Krylov subspace solvers of the DDM linear system~\eqref{ddm} is largely determined by the choice of Robin boundary conditions therein. More effective Robin/impedance boundary conditions on the subdomain interfaces are known to improve the performance of iterative DDM solvers~\cite{steinbach2011stable,boubendirDDM,vion2014double,Gander1}. These generalized Robin boundary conditions consist in replacing the classical $i\eta$ term by operators that approximate the Dirichlet-to-Neumann (DtN) operators of adjacent domains. For instance, the ideal Robin operator on the interface $\Gamma_{12}$ corresponding to the domain $\Omega_1$ consists of the operator $\varepsilon_2^{-1} Y^2|_{\Gamma_{21}}$, where $Y^2$ is the DtN operator corresponding to the domain $\Omega_2$ with zero Dirichlet boundary conditions on $\partial\Omega_2\setminus\Gamma_{21}$.  With this very choice, the ensuing DDM algorithm converges in precisely two iterations~\cite{Nataf}, at least in the case when $\Omega_j$ are half planes.  Similarly, the ideal Robin operator on the interface $\Gamma_{10}$ corresponding to the domain $\Omega_1$ can be shown to consist of the operator $\varepsilon_0^{-1} Y^0|_{\Gamma_{01}}$. Although DtN operators are not always defined for interior subdomains, they are always well defined in the exterior domain $\Omega_0$. However, even when they are properly defined DtN maps are non-local operators whose computation can be expensive. Their computation, whenever possible, can be obtained via boundary integral operators. For instance, using Green's identities in the domain $\Omega_2$ and taking into consideration the null Dirichlet boundary conditions on $\partial\Omega_2\setminus\Gamma_{21}$, we obtain the expression:
\[
u_2=-DL_{\Gamma_{21},2} u_2 + SL_{\partial\Omega_2,2}\partial_{n_2}u_2\quad \mbox{in}\quad\Omega_2,
\]
which  upon application of Dirichlet traces, leads to the identity
\begin{equation}\label{eq:Sinv}
Y^2=S_{\partial\Omega_2,2}^{-1}\left(\frac{1}{2}I+K_{\Gamma_{21},2}\right).
\end{equation}
The invertibility of the operator $S_{\partial\Omega_2,2}$ in the equation above, and hence the well-posedness of the DtN operator $Y^2$, can be guaranteed provided the subdomain $\Omega_2$ is small enough, typically less than one wavelength across. A simple solution that would allow one to consider subdomains of any size is to consider DtN operators $Y^{2,c}$ corresponding to complexified wavenumbers $k_2+i\sigma_2,\ \sigma_2>0$ instead of the operators $Y^2$. Using these operators, we can define a transmission operator on the interface $\partial\Omega_2$ in the form  
\begin{equation}\label{eq:DtNR}
 \mathcal{T}_1^{DtN}=\varepsilon_2^{-1}X_{12}Y^{2,c} + \varepsilon_0^{-1}X_{10}Y^0
\end{equation}
and similar transmission operators on the interfaces $\partial\Omega_2$ and $\partial\Omega_0$ respectively. We then match DtN Robin boundary conditions (DtNR) on the subdomain interfaces in the form 
\begin{equation}\label{eq:DNRBC}
 \varepsilon_1^{-1}\partial_{n_1}u_1+\mathcal{T}_1^{DtN}u_1= -\varepsilon_j^{-1}(\partial_{n_j}u_j+\delta_j^0\partial_{n_j}u^{inc})+\mathcal{T}_1^{DtN}(u_j+\delta_j^0 u^{inc}),\ j\in\{0,2\}.
\end{equation}
\begin{remark}\label{remDtN}
  We point out several issues related to the DDM formulation that uses transmission operators defined in equation~\eqref{eq:DtNR}. The primary issue related to generalized Robin boundary conditions matching is whether these are equivalent to the continuity of the Dirichlet and (scaled) Neumann traces on subdomain interfaces. It is a straightforward matter to see that the aforementioned equivalence is guaranteed as long as the operators $\varepsilon_2^{-1}\Pi_{21}Y^{2,c}+\varepsilon_1^{-1}\Pi_{12}Y^{1,c}$ and their conterparts are invertible in appropriate functional spaces of functions/distributions on $\Gamma_{12}$. These functional spaces are related to the mapping property $\Pi_{21}Y^{2,c}, \Pi_{12}Y^{1,c}:\widetilde{H}^{1/2}(\Gamma_{12})\to H^{-1/2}(\Gamma_{12})$ (we used the same notation and definitions for the Sobolev spaces on the interfaces $\Gamma_{12}$ as in reference~\cite{mclean:2000}). Following similar arguments to those in the proof of Theorem~\ref{thm_inv}, we can establish the invertibility of the operator $\varepsilon_2^{-1}\Pi_{21}Y^{2,c}+\varepsilon_1^{-1}\Pi_{12}Y^{1,c}$. The second issue is related to the mapping property of the operators $\mathcal{T}_1^{DtN}$ defined in equation~\eqref{eq:DtNR} and the well-posedness of the Helmholtz equation in the domain $\Omega_1$ with boundary conditions that feature the generalized Robin operator of the left-hand side of equation~\eqref{eq:DNRBC}. Clearly, the straightforward extension by zero of the operator $\Pi_{21}Y^{2,c}:\widetilde{H}^{1/2}(\Gamma_{12})\to H^{-1/2}(\Gamma_{12})$ (e.g. the operator $X_{12}Y^{2,c}$) does not give rise to an operator that maps $H^{1/2}(\partial\Omega_1)$ to $H^{-1/2}(\partial\Omega_1)$, and thus, the operators $\mathcal{T}_1^{DtN}$ defined in equation~\eqref{eq:DtNR} do not enjoy the latter mapping property. We discuss in what follows a localization/extension procedure via smooth cut-off functions that delivers transmission operator with desired mapping properties.
  \end{remark}
Besides the issues raised in Remark~\ref{remDtN}, the DtN operators for the domains $\Omega_j, 0\leq j$, even when properly defined, are expensive to compute and to apply in equations~\eqref{eq:DNRBC}. We therefore rely on suitable approximations of DtN operators that can be readily computed and at the same time lead to well posed generalized Robin boundary conditions in the domains $\Omega_j,\ 0\leq j$. These approximations start from an exact representation of the DtN operators in each domain. For instance, an application of  Neumann traces to the Green's identities in the domain $\Omega_2$ leads to 
\[
\left(\frac{1}{2}I-K_{\Gamma_{21},2}^\top\right)Y^2 = -N_{\partial\Omega_2,2},
\]
and thus, whenever admissible,
\[
Y^2=-\left(\frac{1}{2}I-K_{\Gamma_{21},2}^\top\right)^{-1}N_{\partial\Omega_2,2}.
\]
An explicit approximation of the DtN operator $Y^2$ can be achieved by keeping the leading order term in the equation above in the sense of pseudodifferential operators
\[
Y^2\approx -2N_{\partial\Omega_2,2}.
\]
The restriction of the operators $N_{\partial\Omega_2,2}$ to the interface $\Gamma_{12}$--which is an open arc--has to be performed with care~\cite{Steinbach,turc2016well}. In addition, the generalized Robin operators should lead to well posed problems in the domain $\Omega_1$. One possibility to blend the approximations of $Y^2|_{\Gamma_{21}}$ and $Y^0|_{\Gamma_{20}}$ that use hypersingular operators is given by~\cite{turc2016well}:
 \begin{equation}\label{eq:calT}
 \mathcal{T}_1=-2\varepsilon_0^{-1}\chi_{01}N_{\partial\Omega_1,k_0+i\sigma_0}\chi_{01}-2\varepsilon_{2}^{-1}\chi_{12}N_{\partial\Omega_1,k_2+i\sigma_2}\chi_{12},\ \sigma_j>0,
 \end{equation}
 where $\chi_{01}$ and $\chi_{12}$ are smooth cutoff functions whose supports lie in $\Gamma_{10}$ and $\Gamma_{12}$, respectively. Given that hypersingular operators are expensive to compute, we can replace the hypersingular operators in equation~\eqref{eq:calT} by principal symbol Fourier multiplier operators. The latter principal symbols are defined as
 \[
 p^N(\xi, k_0+i\sigma_0)=-\frac{1}{2}\sqrt{|\xi|^2-(k_0+i\sigma_0)^2}\quad {\rm and}\quad p^N(\xi, k_2+i\sigma_2)=-\frac{1}{2}\sqrt{|\xi|^2-(k_2+i\sigma_2)^2},
 \]
 where the square root branches are chosen such that the imaginary parts of the principal symbols are positive. The principal symbol Fourier multipliers are defined in the Fourier space $TM(\partial\Omega_1)$~\cite{AntoineX} as
 \begin{equation}
   [PS(N_{\partial\Omega_1,k_j+i\sigma_j})\varphi_1]\hat\ (\xi)=p^N(\xi, k_j+i\sigma_j)\hat{\varphi_1}(\xi)
 \end{equation}
 for a density $\varphi_1$ defined on $\partial\Omega_1$. We define accordingly
 \begin{equation}\label{eq:calPST}
 PS(\mathcal{T}_1)=-2\varepsilon_0^{-1}\chi_{01}PS(N_{\partial\Omega_1,k_0+i\sigma_0})\chi_{01}-2\varepsilon_{2}^{-1}\chi_{12}PS(N_{\partial\Omega_1k_2+i\sigma_2})\chi_{12},\ \sigma_j>0\;,
 \end{equation}
and we use the operators in equation~\eqref{eq:calPST} as generalized Robin/impedance operators. Thus, the classical Robin boundary conditions on $\partial \Omega_1$ are replaced by the following Generalized Square Root Robin boundary conditions (GSqR)
 \begin{equation}\label{eq:GSqR}
 \varepsilon_1^{-1}\partial_{n_1}u_1+PS(\mathcal{T}_1)u_1=-\varepsilon_j^{-1}(\partial_{n_j}u_j+\delta_j^0\partial_{n_j}u^{inc})+PS(\mathcal{T}_1)(u_j+\delta_j^0 u^{inc}),\ j\in\{0,2\}.
 \end{equation}
 We note that the operators $PS(\mathcal{T}_1):H^{1/2}(\partial\Omega_1)\to H^{-1/2}(\partial\Omega_1)$ are coercive, and thus the Helmholtz equations in the domain $\Omega_1$ with generalized Robin boundary conditions defined by the operators in the left-hand side of equation~\eqref{eq:GSqR} are well posed~\cite{turc2016well}. Similar generalized impedance operators can be defined for the domains $\Omega_0$ and $\Omega_2$ and then incorporated in a DDM algorithm that computes the generalized Robin data
 \[
 f_j^g:=(\varepsilon_j^{-1}\partial_{n_j}u_j+PS(\mathcal{T}_j)u_j)|_{\partial\Omega_j},\ 0\leq j
 \]
 by making use of suitably defined generalized RtR maps $\mathcal{S}^{g,j}$.
\begin{remark}\label{rmGDtN}
 While the new blended operators $PS(\mathcal{T}^j),j=0,1,2$ give rise to well-posed Helmholtz equations in each subdomain, the issue of the invertibility of the operator $(PS(\mathcal{T}^1)+PS(\mathcal{T}^2))|_{\Gamma_{12}}$--and hence the equivalence between conditions~\eqref{eq:GSqR} and the continuity of the Dirichlet and Neumann traces on subdomain interfaces--is an open question. In addition, the generalized RtR maps are no longer unitary and thus the proof of solvability of the DDM linear system with GSqR~\eqref{eq:GSqR} also remains open. 
\end{remark}

\section{High-order Nystr\"om discretizations of the MTF equation~\eqref{eq:MTF1}, g-MTF equation~\eqref{eq:sk_system} and the DDM formulation~\eqref{ddm}~\label{Nystrom}}

\parskip 10pt plus2pt minus1pt
\parindent0pt 

We use Nystr\"om discretizations of  the MTF equation~\eqref{eq:MTF1}, g-MTF equation~\eqref{eq:sk_system} and the DDM formulation~\eqref{ddm} that relies on discretizations of the four BIO in the Calder\'on calculus. The latter, in turn, encompass: (a) use of graded meshes based on sigmoid transforms that cluster polynomially discretization points toward corners; (b) splitting of the kernels of the parametrized versions of the boundary integral operators that feature in the  MTF equation~\eqref{eq:MTF1} and the  g-MTF equation~\eqref{eq:sk_system} into sums of regular quantities and products of periodized logarithms and regular quantities; (c) trigonometric interpolation of the densities of the boundary integral operators; and, (d) analytical expressions for the integrals of products of periodic singular and weakly singular kernels and Fourier harmonics. These discretizations were introduced in~\cite{dominguez2015well} where this methodology is presented in full detail. The main idea of our Nystr\"om discretization is to incorporate sigmoid transforms~\cite{KressCorner} in the parametrization of a closed Lipschitz curve $\Gamma$ and then split the kernels of the Helmholtz BIO into smooth and singular components. Using graded meshes that avoid corner points and classical singular quadratures of Kusmaul and Martensen~\cite{kusmaul,martensen}, we employ the Nystr\"om discretization presented in~\cite{dominguez2015well} to produce high-order approximations of the BIO that enter the Calder\'on projectors. We note that the role of weighted traces is to increase the regularity of the densities to be integrated, regularity that can be tuned by simply increasing the polynomial order in the sigmoid transforms.  In what follows, we provide details relevant to each of the three formulations considered.

For each of the subdomains $\Omega_j$, $j=0,1,2$, we employ graded meshes $\mathbf{x}^j_{m},m=0,\ldots,N_j-1$ on $\partial \Omega_j$ with the same polynomial degree of the sigmoid transform. All these meshes are subsequently shifted by the same amount--basically half a discretization step-size--so that none of the grid points corresponds to a triple/multiple junction or a corner point. We note that on a common interface between two subdomains $\Omega_j$ and $\Omega_\ell$ the grid points corresponding to the mesh in each subdomain coincide.

{\em Nystr\"om discretization of the MTF equation~\eqref{eq:MTF1}}. Following the prescriptions in~\cite{dominguez2015well} for the discretization of transmission boundary integral equations for one subdomain, i.e.~no triple/multiple junctions, we define {\em weighted} Neumann traces on the boundary of each subdomain $\Omega_j$
\[
\partial_{n_j}^{w}u_j:=|\mathbf{x}_j'|\partial_{n_j}u_j
\]
where $\mathbf{x}_j$ is a $[0,2\pi]$ parametrization of the boundary $\partial\Omega_j$ that incorporates sigmoid transforms on each of the smooth arcs of $\partial\Omega_j$. This procedure requires introducing {\em weighted} parametrized versions of the four BIO of the Calder\'on calculus; these were discussed in detail in the same reference~\cite{dominguez2015well}. 

{\em Nystr\"om discretization of the g-MTF equation~\eqref{eq:sk_system}}. The additional challenge in the Nystr\"om discretization of the g-MTF consists of the fact that the skeleton $\Gamma$ is no longer a closed curve--recall that the Nystr\"om discretizations that rely on the singular quadratures of Kusmaul and Martensen~\cite{kusmaul,martensen} require that the domain of integration is a closed curve and thus all the integrands to be dealt with are {\em periodic}. Nevertheless, introducing {\em weighted} quantities
\[
v^w:=|\mathbf{x}'| v\qquad p^w:=|\mathbf{x}'|p\qquad {\rm on}\ \Gamma=\Gamma_{01}\cup\Gamma_{12}\cup\Gamma_{02},
\]
we obtain a {\em weighted} g-MTF whose unknowns are $(v^w, p^w)$ by multiplying both sides of equations~\eqref{eq:sk_system} by the weight $|\mathbf{x}'|$. We note that the weight $|\mathbf{x}'|$ is related to a parametrization of the skeleton $\Gamma$ that is easily produced from individual parameterizations of the domains $\partial\Omega_j, j=0,1,2$. We note that the weighted unknowns $(v^w, p^w)$ vanish algebraically at the triple/multiple junction points of the skeleton $\Gamma$. This simple fact allows us to evaluate easily weighted boundary integral operators defined on $\Gamma$. Indeed, the evaluation of a parametrized weighted integral operator of the type
\[
(Iv^w)(\mathbf{x}(t))=\int_\Gamma \mathcal{G}_k^w(\mathbf{x}(t),\mathbf{x}(\tau))v^w(\mathbf{x}(\tau))d\tau,\quad \mathbf{x}(t)\in\Gamma_{01},
\]
can be simply performed as
\begin{equation}\label{eq:split_def}
(Iv^w)(\mathbf{x}(t))=\int_{\partial\Omega_1} \mathcal{G}_k^w(\mathbf{x}(t),\mathbf{x}(\tau))v^w_1(\mathbf{x}(\tau))d\tau + \int_{\partial\Omega_0} \mathcal{G}_k^w(\mathbf{x}(t),\mathbf{x}(\tau))v^w_0(\mathbf{x}(\tau))d\tau,
\end{equation}
where
\[
v^w_1(\mathbf{x}(\tau)):=v^w(\mathbf{x}(\tau)),\ \mathbf{x}(\tau)\in \partial\Omega_1=\Gamma_{01}\cup\Gamma_{12}\qquad v^w_0(\mathbf{x}(\tau)):=\begin{cases}v^w(\mathbf{x}(\tau)),\ &\mathbf{x}(\tau)\in \Gamma_{02} \\ 0 ,\  & \mathbf{x}(\tau)\in \Gamma_{02}\end{cases}.
\]
Each of the two integrals in equation~\eqref{eq:split_def} is performed on a closed curve and involves $2\pi$ periodic and regular functions, and thus can be computed with high-order accuracy by the Nystr\"om method depicted above. All of the weighted integral operators that enter the weighted g-MTF  can be effected in a similar manner.

{\em Discretization of the DDM equation~\eqref{ddm} based on Nystr\"om discretizations of the RtR maps}. The RtR maps $\mathcal{S}^j$ can be computed via well-conditioned BIO and subsequently discretized following the Nystr\"om discretizations introduced in~\cite{turc2016well}. In a nutshell, the method introduced in~\cite{turc2016well} is a direct regularized boundary integral equation formulation of interior/exterior Helmholtz equations with Robin/impedance boundary conditions of the type employed by the DDM formulation~\eqref{ddm}. Specifically, if $u_j$ is the solution of the following Helmholtz problem with Robin boundary conditions, augmented by $u_0$ is radiative in the exterior domain $\Omega_0$,
 \begin{eqnarray*}
   \Delta u_j+k_j^2u_j&=&0\ {\rm in}\ \Omega_j\\
   \varepsilon_j^{-1}\partial_{n_j} u_j+i\eta u_j&=&\psi_j\ {\rm on}\ \partial \Omega_j,
 \end{eqnarray*}
then $u_j|_{\partial\Omega_j}$ can be shown to satisfy the following well-conditioned regularized boundary integral equation~\cite{turc2016well}:
\begin{eqnarray}\label{eq:CFIER2}
\mathcal{A}_{j}(u_j|_{\partial\Omega_j})&=&\varepsilon_j(S_j+S_{\kappa_j}-2S_{\kappa_j} K_j^\top)\psi_j,\quad \kappa_j=k_j+i\sigma_j,\ \sigma_j>0,\nonumber\\
\mathcal{A}_{j}&:=&\frac{1}{2}I-2S_{\kappa_j} N_j+ S_{\kappa_j}Z_j-2 S_{\kappa_j} K_j^\top Z_j +K_j+S_jZ_j, 
\end{eqnarray}
where $Z_j=i\eta\varepsilon_j$. Thus the RtR operator $\mathcal{S}^j$ can be expressed as
\begin{equation}\label{eq:SjBI}
  \mathcal{S}^j=I-2Z_j\mathcal{A}_j^{-1}(S_j+S_{\kappa_j}-2S_{\kappa_j} K_j^\top).
  \end{equation}
The generalized Robin conditions~\eqref{eq:DtNR} and~\eqref{eq:GSqR} are amenable to the same treatment provided in~\cite{turc2016well}, the only difference being that the scalar impedance $Z_j$ in equations~\eqref{eq:CFIER2} and~\eqref{eq:SjBI} is replaced by impedance operators defined in equations~\eqref{eq:DNRBC} and~\eqref{eq:calPST} respectively. In order to avoid complications related to singularities at junction/cross points, we replace in the DDM algorithm the RtR maps by {\em weighted} parametrized counterparts
\[
\mathcal{S}^{j,w}(\varepsilon_j^{-1}|\mathbf{x}_j'|\partial_{n_j}u_j-i\eta\ u_j):=\varepsilon_j^{-1}|\mathbf{x}_j'|\partial_{n_j}u_j+i\eta\ u_j.
\]
Collocated discretizations of the latter weighted RtR maps can be easily computed through a simple modification of the methodology introduced in~\cite{turc2016well} and recounted above. However, the representation of RtR maps in terms of BIO requires use of inverses of the operators $\mathcal{A}_j$ cf. equation~\eqref{eq:SjBI}. In order for the DDM algorithm to be efficient, the electromagnetic/acoustic size of subdomains $\Omega_j$ should be amenable to application of direct linear algebra solvers for calculations of the inverses of the collocation of the matrices $\mathcal{A}_j$. The discretization of the weighted RtR maps corresponding to each domain $\partial \Omega_j$ is thus constructed as $N_j\times N_j$ collocation matrices  $\mathcal{S}^j_N$. Alternatively, the application of the collocated versions of the operators $\mathcal{S}^j$ can be effected via Krylov subspace iterative solutions of equations~\eqref{eq:CFIER2} for each domain $\Omega_j$. However, the numbers of iterations needed for solution of equations~\eqref{eq:CFIER2} does grow with the electromagnetic/acoustic size of the interior bounded domains $\Omega_j$, $j\neq0$~\cite{turc2016well}. Therefore, it is advisable to further subdivide the interior domains $\Omega_j$ when their size is deemed too large for efficient computations of their corresponding RtR maps. In the case of the exterior/unbounded domain $\Omega_0$, the numbers of iterations required for the solution of equations~\eqref{eq:CFIER2} are small regardless of the frequency $\omega$~\cite{turc2016well}, and thus the application of $\mathcal{S}^0_N$ can be performed effciently by iterative methods.

Once the collocation matrices $\mathcal{S}^j_N$ are constructed, the discretization of the interface subdomain RtR maps (e.g.~the maps $\mathcal{S}^j_{j\ell,\ell j}$ and all the other ones defined in equation~\eqref{SL}) is straightforward since it simply amounts to extracting suitable blocks from the matrices $\mathcal{S}^j_N$. Indeed, the discretization of the operators $\mathcal{S}^j_{j\ell,\ell j}$ consists of extracting from the collocation matrix $\mathcal{S}^j_N$ the block that corresponds to self-interactions of the grid points on the interface/edge $\partial \Omega_j\cap\partial \Omega_\ell$. This is possible since none of these mesh points $\mathbf{x}^j_m$ corresponds to a corner of $B_j$. We also mention that all the matrix inversions needed in the Schur complement algorithm (\emph{cf.}~formula~\eqref{inv_matrix_explicit}) are performed by either direct linear algebra methods: the size of these matrices does not exceed the number of unknowns need to discretize $\partial \Omega_0$ or by Neumann series (\emph{cf.} equation~\eqref{eq:Neumannseries}). In case when the large domains $\Omega_j$, $j\neq 0$, are further subdived into smaller non-overlapping subdomains $\Omega_j=\cup_{n=1}^J\Omega_{j_n}$ for which RtR maps $\mathcal{S}^{j_n}_{N_{j_n}}$ can be computed via direct linear algebra methods. The RtR maps $\mathcal{S}^j_N$ themselves can be computed by merging the maps $\mathcal{S}^{j_n}_{N_{j_n}}$. However, as pointed in Table~\ref{comp5} in Section~\ref{num}, the numbers of iterations needed by an iterative solver of the classical DDM system~\eqref{ddm} as well as the DtNR and GSqR DDM system do grow with the number of additional subdomains.  

The comments above are also applicable to the DtNR and GSqR DDM algorithms: the size of the interior subdomains has to be amenable to direct solvers computations of the generalized RtR maps. The DtNR DDM algorithm require that the subdomain DtN maps are precomputed, which is also done via Nystr\"om solvers of the boundary integral formulation~\eqref{eq:Sinv}. The discretization of the operators $\mathcal{T}_j^{DtN}$defined in equations~\eqref{eq:DtNR} is obtained by simply extracting suitable blocks from the collocation matrices corresponding to discretized DtN maps. These blocks are subsequently used in the computation of collocated generalized RtR  maps. We note that given that the Fourier multiplier operators defined in equations~\eqref{eq:calPST} can be efficiently discretized via FFTs~\cite{dominguez2015well}, the cost of computing subdomain generalized RtR maps is virtually identical to that of computing subdomain classical RtR map. This is in contrast to volumetric solvers for local domain problems (e.g.~finite difference/finite element) where Pad\'e approximations need to be employed~\cite{vion2014double} for approximating the square root operators in equations~\eqref{eq:calPST}. Both DtNR and GSqR DDM algorithms require significantly fewer Krylov subspace iterations than the classical DDM algorithm, and the reduction in number of iterations is more striking in the high-contrast cases.   

 \section{Numerical results\label{num}}

 We present in what follows numerical experiments with various formulations of the transmission problems~\eqref{system_tj} considered in this text. We use as test configurations with piece-wise constant dielectric permittivities $\varepsilon_j$ in geometries  depicted in Figures~\ref{fig:tj}, in which the circular domain has radius one. All of the formulations considered were discretized following the prescriprion in Section~\ref{Nystrom}. We used discretizations that place the same number of discretization points on each interface--matching boundary meshes, and sigmoid transforms of polynomial degree 3 for the MTF equation~\eqref{eq:MTF1} and the DDM formulation~\eqref{ddm} with $\eta=k_0$, and respectively degree 4 for the g-MTF~\eqref{eq:sk_system}. We chose $\sigma_j=0.1$ in the definition of the operators~\eqref{eq:DtNR} incorporated in the DtNR and, respectively, $\sigma_j=k_j^{1/3}$ in the definition of the operators~\eqref{eq:calPST} that are incorporated in GSqR. We mention that the more robust choice of operators in equation~\eqref{eq:DtNR} gives rise to DDM whose iterative performance is identical to that of DDM that use exact DtN maps. In all the numerical experiments we present numbers of GMRES iterations for various solvers to reach a residual of $10^{-4}$ and present errors in the far-field for $1,024$ equi-spaced far-field directions. The far-field errors were computed using fine discretizations of the MTF solver. As previously mentioned, the discretizations of both g-MTF and the DDM formulation use the same number of unknowns, whereas that of MTF uses twice as many unknowns.  

On account of MTF formulation~\eqref{eq:MTF1} being an integral formulation of the first kind, the number of iterations needed by Krylov subspace iterative solvers (e.g.~GMRES) of its discretization can be considerable, especially for high-contrast configurations at high-frequencies. Fortunately, simple and quite effective preconditioners of the MTF formulations are readily available. Indeed, a simple preconditioner of the MTF fomulation~\eqref{eq:MTF1} consists of the Calder\'on diagonal preconditioner~\cite{jerez-hanckes1}:
\begin{equation}\label{eq:MTF1C}
  \begin{bmatrix}\mathcal{C}_0 & 0 & 0 \\ 0 & \mathcal{C}_1 & 0 \\ 0 & 0 & \mathcal{C}_2\end{bmatrix}\begin{bmatrix}
    \mathcal{C}_0& \mathbb{X}_{01} & \mathbb{X}_{02}\\
    \mathbb{X}_{10} & \mathcal{C}_1 & \mathbb{X}_{12}\\
    \mathbb{X}_{20} & \mathbb{X}_{21} & \mathcal{C}_2
  \end{bmatrix}
    \begin{bmatrix}\mathbf{u}_0\\ \mathbf{u}_1\\ \mathbf{u}_2\end{bmatrix}=\frac{1}{2}\begin{bmatrix}\mathcal{C}_0 & 0 & 0 \\ 0 & \mathcal{C}_1 & 0 \\ 0 & 0 & \mathcal{C}_3\end{bmatrix} \begin{bmatrix}\mathbf{u}^{inc}\\ -\mathbb{X}_{10}\mathbf{u}^{inc} \\ -\mathbb{X}_{20}\mathbf{u}^{inc}\end{bmatrix}.
\end{equation}
Another possibility is to eliminate the unknown $\mathbf{u}_0$ from the MTF formulation~\eqref{eq:MTF1} using Schur complements. Indeed, using the fact that $\mathcal{C}_0$ is a projector, we get that $\mathbf{u}_0$ can be explictly eliminated from the MTF system~\eqref{eq:MTF1} via the equation 
\[
\mathbf{u}_0=2\mathcal{C}_0\mathbf{u}^{inc}-4\mathcal{C}_0\mathbb{X}_{01}\mathbf{u}_1-4\mathcal{C}_0\mathbb{X}_{02}\mathbf{u}_2,
\]
leading to a reduced MTF formulation:
\begin{equation}\label{eq:MTF1S}
  \begin{bmatrix}\mathcal{C}_{1}-4\mathbb{X}_{10}\mathcal{C}_0\mathbb{X}_{01} & \mathbb{X}_{12}-4\mathbb{X}_{10}\mathcal{C}_0\mathbb{X}_{02}\\
\mathbb{X}_{21} - 4\mathbb{X}_{20}\mathcal{C}_0\mathbb{X}_{01} & \mathcal{C}_2-4\mathbb{X}_{20}\mathcal{C}_0\mathbb{X}_{02}
\end{bmatrix} \begin{bmatrix}\mathbf{u}_1\\ \mathbf{u}_2\end{bmatrix}=-\frac{1}{2}\begin{bmatrix}\mathbb{X}_{10}\mathbf{u}^{inc} \\ \mathbb{X}_{20}\mathbf{u}^{inc}\end{bmatrix}-2\begin{bmatrix}\mathbb{X}_{10}\mathcal{C}_0\mathbf{u}^{inc} \\ \mathbb{X}_{20}\mathcal{C}_0\mathbf{u}^{inc}\end{bmatrix}.
\end{equation}
Formulation~\eqref{eq:MTF1S}, in turn, can be further preconditioned by the operator $diag(\mathcal{C}_1,\mathcal{C}_2)$, a formulation for which we report the numbers of GMRES iterations in this section.

We start the presentation of numerical results by illustrations in Tables~\ref{comp1} and~\ref{comp2} of the performance of the three types of solvers dicussed in this text in the case of multiple junction configurations with very high-contrast materials at low incident frequencies $\omega$. We continued in Tables~\ref{comp3} and~\ref{comp4} with results concerning multiple junction configurations at high frequencies. Finally, we conclude in Table~\ref{comp5} with an illustration of the behavior of the Nystr\"om solvers of the formulations considered as the size of the discretization is increased. The accuracy achieved by the GSqR DDM solver is at the same level to that of the DDM solver in all of the cases considered, and for this reason and in order to better streamline the information presented we chose not to report it. 

 \begin{table}
   \begin{center}
     \resizebox{!}{0.8cm}
{
\begin{tabular}{|c|c|c|c|c|c|c|c|c|c|c|c|c|c|c|}
\hline
$\omega$ & \multicolumn{4}{c|} {DDM~\eqref{DDM_explicit}} & \multicolumn{2}{c|} {DDM Schur~\eqref{eq:final}} & \multicolumn{2}{c|} {DDM Schur~\eqref{eq:Neumannseries}} & \multicolumn{4}{c|} {MTF}& \multicolumn{2}{c|} {g-MTF~\eqref{eq:sk_system}}\\
\cline{2-15}
& It & It DtNR & It GSqR & $\varepsilon_\infty$ & It & $\varepsilon_\infty$ & It & $\varepsilon_\infty$ & It~\eqref{eq:MTF1} & It~\eqref{eq:MTF1C} & It~\eqref{eq:MTF1S}& $\varepsilon_\infty$ &It&$\varepsilon_\infty$\\
\hline
1 & 137 & 16 & 25 & 4.8 $\times$ $10^{-3}$ & 15 & 4.9 $\times$ $10^{-3}$ & 15 & 3.2 $\times$ $10^{-2}$ & 431 & 242 & 160 & 5.3 $\times$ $10^{-3}$&159&3.9$\times 10^{-3} $\\
2 & 182 & 19 & 40 & 1.8 $\times$ $10^{-3}$ & 18 & 2.6 $\times$ $10^{-3}$ & 19 & 4.8 $\times$ $10^{-2}$ &750 & 433 & 293 & 2.4 $\times$ $10^{-3}$&270&3.2$\times 10^{-3}$ \\
4 & 282 & 19 & 68 & 2.9 $\times$ $10^{-3}$ & 20 & 1.3 $\times$ $10^{-3}$ & 21 & 1.0 $\times$ $10^{-2}$ & 1,335 & 732 & 508 & 3.4 $\times$ $10^{-3}$&430&2.8$\times 10^{-3}$ \\
\hline
\end{tabular}
}
\caption{Performance of the various formulations considered in this text in the three subdomain case with $\varepsilon_0=1$, $\varepsilon_1=64$, and $\varepsilon_2=256$. The numbers of unknowns required by the DDM and g-MTF are $384$, $768$, and $1,536$ respectively; the MTF uses twice as many unknowns in each case. The number of terms used in the Neumannn series~\eqref{eq:Neumannseries} were $N_{trunc}=40$, $N_{trunc}=80$, and $N_{trunc}=80$ in the cases $\omega=1$, $\omega=2$, and $\omega=4$ respectively.\label{comp1}}
\end{center}
 \end{table}

 \begin{table}
   \begin{center}
     \resizebox{!}{0.6cm}
{
\begin{tabular}{|c|c|c|c|c|c|c|c|c|c|c|c|c|c|c|}
\hline
$\omega$ & \multicolumn{4}{c|} {DDM~\eqref{DDM_explicit}} & \multicolumn{2}{c|} {DDM Schur~\eqref{eq:final}} & \multicolumn{2}{c|} {DDM Schur~\eqref{eq:Neumannseries}} & \multicolumn{4}{c|} {MTF}&\multicolumn{2}{c|}{g-MTF~\eqref{eq:sk_system}}\\
\cline{2-15}
& It & It DtNR & It GSqR & $\varepsilon_\infty$ & It & $\varepsilon_\infty$ & It~\eqref{eq:MTF1} & It~\eqref{eq:MTF1C} & It~\eqref{eq:MTF1S}& $\varepsilon_\infty$ & It&$\varepsilon_\infty$& It&$\varepsilon_\infty$\\
\hline
1 & 977 & 39 & 70 & 1.6 $\times$ $10^{-3}$ & 12 & 1.9 $\times$ $10^{-3}$ & 14 & 3.5 $\times$ $10^{-2}$ &2,419 & 1,296 & 1,107 & 4.5 $\times$ $10^{-3}$ &942&1.1$\times$ $10^{-3}$\\
2 & 1,788 & 72 & 140 & 2.9 $\times$ $10^{-3}$ & 20 & 3.3 $\times$ $10^{-3}$ & 23 & 2.8 $\times$ $10^{-2}$ & 4,534 & 2,442 & 2,154 & 5.4 $\times$ $10^{-3}$&1,782&2.1$\times$ $10^{-3}$ \\
\hline
\end{tabular}
}
\caption{Performance of the various formulations considered in this text in the five subdomain case with $\varepsilon_0=1$, $\varepsilon_1=64$, $\varepsilon_2=256$, $\varepsilon_3=1024$, and $\varepsilon_4=4096$. The numbers of unknowns required by the DDM and g-MTF are $4,608$ and $9,216$ respectively; the MTF uses twice as many unknowns in each case. The number of terms used in the Neumannn series~\eqref{eq:Neumannseries} were $N_{trunc}=160$ in both cases $\omega=1$,  and $\omega=2$ respectively.\label{comp2}}
\end{center}
 \end{table}

 \begin{table}
   \begin{center}
     \resizebox{!}{0.8cm}
{
\begin{tabular}{|c|c|c|c|c|c|c|c|c|c|c|c|c|c|c|}
\hline
$\omega$ & \multicolumn{4}{c|} {DDM~\eqref{DDM_explicit}} & \multicolumn{2}{c|} {DDM Schur~\eqref{eq:final}} & \multicolumn{2}{c|} {DDM Schur~\eqref{eq:Neumannseries}} & \multicolumn{4}{c|} {MTF}&\multicolumn{2}{c|} {g-MTF~\eqref{eq:sk_system}}\\
\cline{2-15}
& It & It DtNR & It GSqR & $\varepsilon_\infty$ & It & $\varepsilon_\infty$ & It & $\varepsilon_\infty$ & It~\eqref{eq:MTF1} & It~\eqref{eq:MTF1C} & It~\eqref{eq:MTF1S}& $\varepsilon_\infty$&It&$\varepsilon_\infty$\\
\hline
4 & 88 & 25 & 27 & 3.6 $\times$ $10^{-3}$ & 25 & 3.5 $\times$ $10^{-3}$ & 26 & 8.5 $\times$ $10^{-3}$ & 212 & 150 & 109 & 4.3 $\times$ $10^{-3}$&108&1.8$\times 10^{-3}$ \\
8 & 142 & 40 & 41 &4.4 $\times$ $10^{-3}$ & 40 & 4.3 $\times$ $10^{-3}$ & 42 & 4.2 $\times$ $10^{-2}$ & 342 & 235 & 176 & 5.8 $\times$ $10^{-3}$&174&4.1$\times 10^{-3}$ \\
16 & 176 & 49 & 57 & 3.1 $\times$ $10^{-3}$ & 45 & 3.0 $\times$ $10^{-3}$ & 48 & 4.9 $\times$ $10^{-2}$ & 483 & 361 & 267 & 6.3 $\times$ $10^{-3}$&280&5.8$\times 10^{-3}$ \\
32 & 274 & 64 & 91 &1.9 $\times$ $10^{-3}$ & 59 & 1.3 $\times$ $10^{-3}$ & 64 & 4.5 $\times$ $10^{-2}$ & 808 & 565 & 405 & 3.8 $\times$ $10^{-3}$ &475&7.8$\times 10^{-3}$\\
\hline
\end{tabular}
}
\caption{Performance of the various formulations considered in this text in the three subdomain case with $\varepsilon_0=1$, $\varepsilon_1=4$, and $\varepsilon_2=16$. The numbers of unknowns required by the DDM and g-MTF are $384$, $768$, $1,536$, and $3,072$ respectively; the MTF uses twice as many unknowns in each case. The number of terms used in the Neumannn series~\eqref{eq:Neumannseries} were $N_{trunc}=20$, $N_{trunc}=20$, $N_{trunc}=30$, and $N_{trunc}=80$ in the cases $\omega=4$, $\omega=8$, $\omega=16$, and $\omega=32$ respectively. \label{comp3}}
\end{center}
 \end{table}

 \begin{table}
   \begin{center}
     \resizebox{!}{0.8cm}
{
\begin{tabular}{|c|c|c|c|c|c|c|c|c|c|c|c|c|c|c|}
\hline
$\omega$ & \multicolumn{4}{c|} {DDM~\eqref{DDM_explicit}} & \multicolumn{2}{c|} {DDM Schur~\eqref{eq:final}} & \multicolumn{2}{c|} {DDM Schur~\eqref{eq:Neumannseries}} & \multicolumn{4}{c|} {MTF}&\multicolumn{2}{c|} {g-MTF~\eqref{eq:sk_system}}\\
\cline{2-15}
& It & It DtNR & It GSqR & $\varepsilon_\infty$ & It & $\varepsilon_\infty$ & It & $\varepsilon_\infty$ & It~\eqref{eq:MTF1} & It~\eqref{eq:MTF1C} & It~\eqref{eq:MTF1S}& $\varepsilon_\infty$ &It&$\varepsilon_\infty$\\
\hline
4 & 427 & 52 & 67 & 2.3 $\times$ $10^{-3}$ & 23 & 2.8 $\times$ $10^{-3}$ & 30 & 3.6 $\times$ $10^{-2}$ & 1,152 & 673 & 599 & 6.0 $\times$ $10^{-3}$ &859&5.4$\times 10^{-3}$\\
8 & 613 & 80 & 109 & 1.1 $\times$ $10^{-3}$ & 32 & 1.3 $\times$ $10^{-3}$ & 40 & 4.5 $\times$ $10^{-2}$ & 2,058 & 1,227 & 1,111 & 6.7 $\times$ $10^{-3}$&1,426&4.9$\times 10^{-3}$ \\
16 & 784 & 106 & 178 & 9.3 $\times$ $10^{-4}$ & 41 & 2.1 $\times$ $10^{-3}$ & 51 & 3.1 $\times$ $10^{-2}$ & 3,371 & 2,163 & 1,997 & 5.3 $\times$ $10^{-3}$&2,471&5.4$\times 10^{-3}$ \\
\hline
\end{tabular}
}
\caption{Performance of the various formulations considered in this text in the five subdomain case with  $\varepsilon_0=1$, $\varepsilon_1=4$, $\varepsilon_2=16$, $\varepsilon_3= 64$, and $\varepsilon_4=256$. The numbers of unknowns required by the DDM and g-MTF are $1,152$, $2,304$ and $4,608$ respectively; the MTF uses twice as many unknowns in each case. The number of terms used in the Neumannn series~\eqref{eq:Neumannseries} were $N_{trunc}=80$, $N_{trunc}=80$, and $N_{trunc}=160$ in the cases $\omega=4$, $\omega=8$, and $\omega=16$ respectively.\label{comp4}}
\end{center}
 \end{table}

 \begin{table}
   \begin{center}
  \resizebox{!}{0.7cm}
{   
\begin{tabular}{|c|c|c|c|c|c|c|c|c|c|c|c|}
\hline
\multicolumn{2}{c|} {DDM~\eqref{DDM_explicit}} & {DDM Schur~\eqref{eq:final}} & \multicolumn{2}{c|} {4 DDM~\eqref{DDM_explicit}} & {4 DDM Schur~\eqref{eq:final}} & \multicolumn{4}{c|} {MTF} & \multicolumn{2}{c|} {g-MTF~\eqref{eq:sk_system}}\\
\hline
It & $\varepsilon_\infty$ & It & It & $\varepsilon_\infty$ & It & It~\eqref{eq:MTF1} & It~\eqref{eq:MTF1C} & It~\eqref{eq:MTF1S}& $\varepsilon_\infty$ &It&$\varepsilon_\infty$\\
\hline
77 & 4.1 $\times$ $10^{-3}$ & 22 & 107 & 4.5 $\times$ $10^{-3}$ & 23 & 142 & 102 & 77 & 3.6 $\times$ $10^{-3}$& 67 &1.4$\times 10^{-3}$ \\
73 & 3.4 $\times$ $10^{-4}$ & 22 & 118 & 3.5 $\times$ $10^{-4}$ & 22 & 145 & 104 & 78 & 2.8 $\times$ $10^{-4}$& 69 &1.7$\times 10^{-4}$ \\
\hline
\end{tabular}
}
\caption{Performance of the various formulations considered in this text in the three subdomain case with $\varepsilon_0=1$, $\varepsilon_1=4$, $\varepsilon_2=16$ and $\omega=2$. The numbers of unknowns required by the DDM and g-MTF are $192$ and $384$ respectively; the MTF uses twice as many unknowns in each case. We present results for the DDM formulations in the case when the subdomains $\Omega_1$ and $\Omega_2$ are each subdivided into two subdomains, and we refer to these cases by 4 DDM.\label{comp5}}
\end{center}
 \end{table}

 \begin{figure}
\centering
\includegraphics[height=50mm]{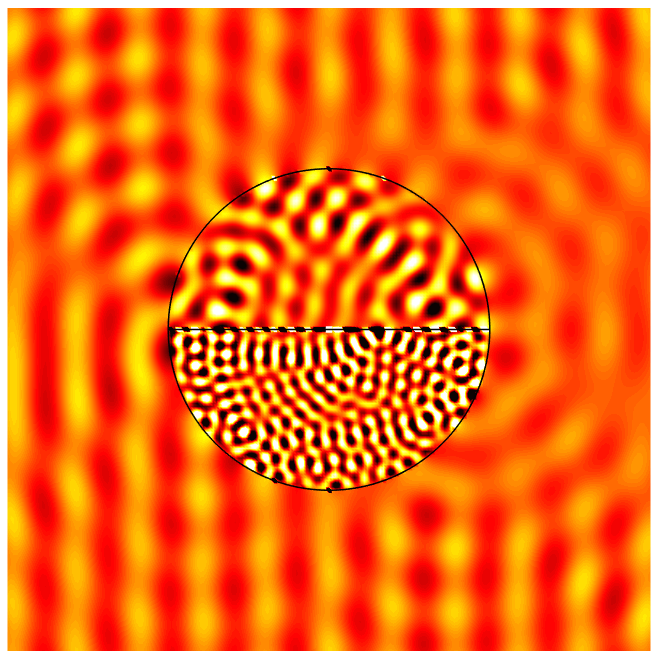}\includegraphics[height=50mm]{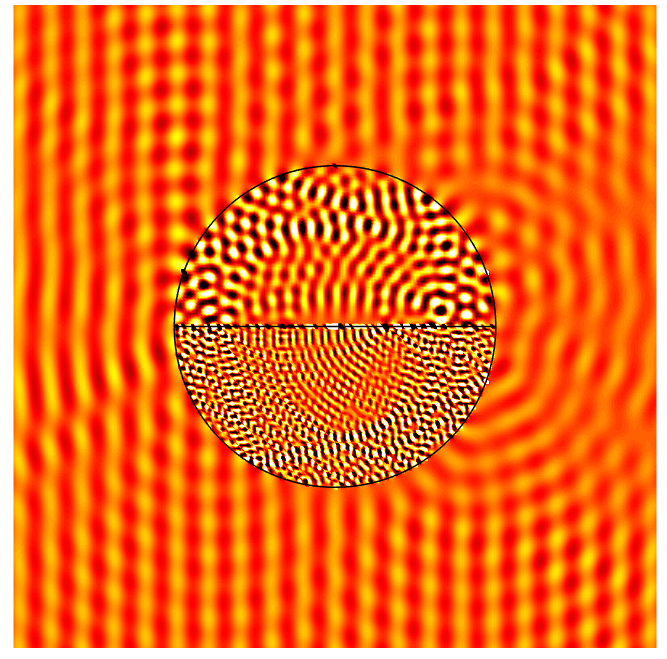}\\
\includegraphics[height=50mm]{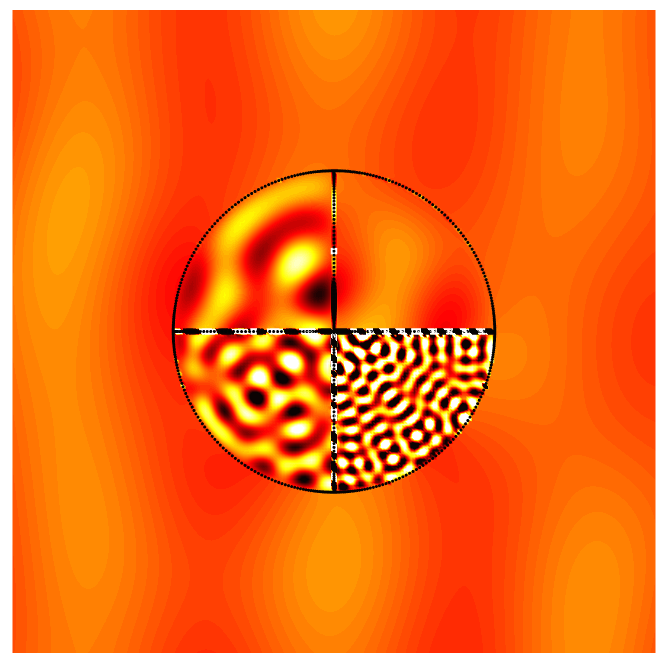}\includegraphics[height=50mm]{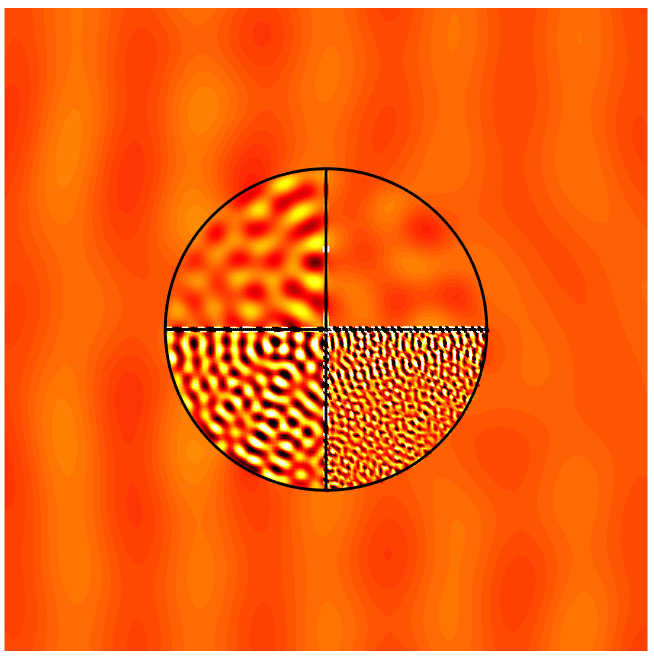}\\
\caption{Plots of the near fields scattered by multiple junction configurations under plane wave incidence of direction $d=(1,0)$. Top panel: $\varepsilon_0=1$, $\varepsilon_1=4$, and $\varepsilon_2=16$; $\omega=16$ (left) and $\omega=32$ (right). Bottom panel: $\varepsilon_0=1$, $\varepsilon_1=4$, $\varepsilon_2=16$, $\varepsilon_3=64$, and $\varepsilon_4=256$; $\omega=4$ (left) and $\omega=8$ (right)}
\label{fig:plots}
 \end{figure}

 {\em Discussion of the performance of various formulations}. We begin by noting that the computational cost of setting up the MTF equation~\eqref{eq:MTF1} and the g-MTF~\eqref{eq:sk_system} is comparable, as the building blocks of these formulations are the same BIO defined on the boundaries of the subdomains. The behavior of iterative solvers for the solution of Nystr\"om discretizations of the two formulations is quite different, given that MTF is a first kind formulation that requires twice as many unknowns than the second kind formulation CFIESK. Calder\'on preconditioning strategies that lead to the MTF formulation~\eqref{eq:MTF1C} and Schur complement strategies combined with  Calder\'on preconditioning strategies that lead to the MTF formulation~\eqref{eq:MTF1S} are effective in reducing the numbers of GMRES iterations and do not require significant additional computational costs. In contrast, the g-MTF~\eqref{eq:sk_system}, whose motivation was canceling the adverse effects of hypersingular operators, is not easily amenable to preconditioning. However, even with preconditioning, both formulations require large numbers of GMRES iterations, especially in the case of high-contrast multi-domain configurations at high-frequencies. 

 As it can be seen from the results presented in this section, solvers based on DDM formulations~\eqref{ddm} and the DDM formulation that incorporate either the DtNR or GSqR give rise to smaller numbers of GMRES iterations than those based on the MTF and g-MTF. We note that effecting the matrix-vector product associated with the DDM equation~\eqref{ddm}, requires in turn application of each of the RtR maps $\mathcal{S}^j$. The latter can be done iteratively via solutions of Robin/impedance boundary values problems in each corresponding domain $\Omega_j$ via equations~\eqref{eq:CFIER2}. We present in Table~\ref{comp6} the numbers of GMRES iterations required for the solution of equations~\eqref{eq:CFIER2} in each subdomain $\Omega_j$ when using a residual of $10^{-4}$. We observe that these numbers grow with the frequency $\omega$, especially in the case of the interior domains $\Omega_1$ and $\Omega_2$. Further subdivisions of the interior domains bring down these numbers, but only modestly so, while increasing considerably the numbers of GMRES iterations needed for the solution of the global DDM system~\eqref{ddm}. 
 
 These findings corroborate those reported in the DDM literature~\cite{steinbach2011stable} and suggest that DDM algorithm are efficient when the size of subdomains is amenable to use of direct solvers for computation of the action of subdomains RtR maps. The Nystr\"om discretizations presented in Section~\ref{Nystrom} constitute very efficient means to discretize both classical and generalized RtR operators. Our DDM algorithms thus proceed in two stages: (1) a precomputation/offline stage where all the subdomain RtR maps are precomputed, and (2) an online stage where the DDM linear system is solved iteratively in a matrix-free fashion. We note that the use of GSqR~\eqref{eq:calPST} in the DDM algorithm gives rise to important savings in numbers of iterations with virtually no additional cost. The numbers of iterations is further decreased in the DtNR DDM algorithm, yet at the expense of adding one more offline stage (0) where the subdomain DtN maps are precomputed. However, the number of subdomains adversely affects the rate of convergence of iterative solvers of the DDM system~\eqref{ddm}, regardless of the type of Robin boundary conditions used.  These findings constitute somewhat of a surprise especially in the case when DtNR boundary conditions~\eqref{eq:DNRBC} are used: one would expect that the incorporation of (almost) exact DtN should result in numbers of DDM iterations that are insensitive to the frequency and the number of subdomains. A similar situation was analyzed in the one dimensional case with constant wavenumber in~\cite{vion2014double}). The use of double sweep preconditioners in slab-type DDM was shown in~\cite{vion2014double} to scale with the frequency and the number of subdomains in the absence of sharp interface sof material discontinuity. The incorporation of double sweep preconditioners for the GSqR DDM is currently underway.  
 
\begin{table}
  \begin{center}
     \resizebox{!}{1.4cm}
{ 
\begin{tabular}{|c|c|c|c|c|c|c|c|c|c|c|c|c|}
\hline
$\omega$ & \multicolumn{6}{c|} {DDM~\eqref{DDM_explicit}}&  \multicolumn{6}{c|} {4 DDM~\eqref{DDM_explicit}}\\
\cline{2-11}
\hline
& It $\Omega_0$ & It $\Omega_1$ & It $\Omega_2$ & It & It DtNR & It GSqR & It $\Omega_0$ & It $\Omega_{11}\ \Omega_{12}$ & It $\Omega_{21}\  \Omega_{22}$ & It & It DtNR & It GSqR \\
\hline
1 & 8 & 30 & 48 & 137 & 16 & 25 & 8 & 27 & 39 & 225 & 42 & 44 \\
2 & 9 & 47 & 79 & 182 & 19 & 40 & 9 & 39 & 64 & 361 & 61 & 73\\
4 & 11 & 70 & 121 & 282 & 19 & 68 & 10 & 59 & 109 & 597 & 104 & 129\\
\hline
4 & 10 & 21 & 36 & 88 & 31 & 30& 10 & 19 & 32 & 153 & 46 & 37 \\
8 & 13 & 25 & 46 & 142 & 40 & 54 & 12 & 21 & 40 & 231 & 64 & 63 \\
16 & 16 & 32 & 60 & 176 & 49 & 77 & 15 & 31 & 56 & 257 & 78 & 78\\
32 & 19 & 39 & 73 & 274 & 64 & 106 & 18 & 38 & 70 & 367 & 101 & 130\\
\hline
\end{tabular}
}
\caption{Left panel: numbers of GMRES iterations needed to reach a residual of $10^{-4}$ for the solution of equation~\eqref{eq:CFIER2} in each subdomain and respectively for the solution of the global DDM system~\eqref{ddm} in the three subdomain case with $\varepsilon_0=1$, $\varepsilon_1=64$, and $\varepsilon_2=256$ (top three rows) and $\varepsilon_0=1$, $\varepsilon_1=4$, and $\varepsilon_2=16$ (bottom four rows). Right panel: same numbers of GMRES iterations when the subdomains $\Omega_1$ and $\Omega_2$ are each subdivided into two half-sized non-overlapping subdomains $\Omega_1=\Omega_{11}\cup \Omega_{12}$ and $\Omega_2=\Omega_{21}\cup\Omega_{22}$ respectively.\label{comp6}}
\end{center}
 \end{table}

The Schur complement DDM algorithm does require that the collocation matrices $\mathcal{S}^j_N$ be precomputed for all the interior subdomains, i.e.~$j\neq 0$. This can be done by either inverting matrices of sizes $N_j\times N_j$, $j\neq 0$ at a cost $N_j^3$, or, alternatively, by further subdividing the domain $\Omega_j$ into smaller subdomains --for which the corresponding collocated RtR maps can be computed effciently with direct methods-- and then applying the merging procedure to the RtR maps of the smaller subdomains. The latter algorithm can also be computationally more advatangeous than the direct approach, especially in the case of identical subdomains. The Schur complement DDM algorithm requires at each stage inversion of matrices (\emph{cf.}~formula~\eqref{inv_matrix_explicit}) whose size depends on the number of unknowns on the common interfaces between subdomains. The matrix inversions required by the merging procedure can be done using direct linear algebra method or using Neumann series~\eqref{eq:Neumannseries}.  If the interior domain $\cup_{1\leq j}\Omega_j$ is decomposed into a collection $P\times P$ subdomains, and $n$ discretization points are used per each subdomain interface, then the cost of hierarchically merging of the subdomain RtR maps to deliver the interior RtR map $\mathcal{S}^{int}$ is $\mathcal{O}(P^3n^3)$. In contrast, direct solvers for MTF would require computational costs of the order $\mathcal{O}(2^9P^6n^3)$, whereas  direct solvers for g-MTF would require costs of the order $\mathcal{O}(2^6P^6n^3)$. If the aforementioned matrix inversions are bypassed via truncated Neumann series, the number of terms needed in the Neumann series is quite small in the case of low contrast/low frequencies, but does grow in the case of very high-contrast/ high frequencies. 

Finally, after the merging of the RtR maps of interior domains $\mathcal{S}^{int}$ has been computed via hierarchical Schur complements, the equation~\eqref{eq:final} is solved either via direct linear algebra solvers or iteratively. If the latter approach is used, the exterior RtR map $\mathcal{S}^0$ needs not be precomputed. Furthermore, the computation of a matrix vector product associated with the exterior RtR map $\mathcal{S}_0$ can be performed efficiently via iterative solutions of equation~\eqref{eq:CFIER2} according to the results in Table~\ref{comp6}. As it can be seen from the results presented in this section, the number of GMRES iterations needed to solve equation~\eqref{eq:final} is quite reasonable throughout the range of problems considered, and does not appear to grow dramatically with the frequency. Nevertheless, equation~\eqref{eq:final} shares the same features with integral equations of the first kind, alternative robust interior/exterior coupling strategies are desirable. One possibility currently under investigation is to modify the generalized Robin operators~\eqref{eq:calPST} so that the same complex wavenumber is used on all the interior interfaces. We conclude by mentioning that, just like the classical iterative DDM algorithm, the Schur complement DDM algorithm is highly parallelizable.

\section{Conclusions}\label{conclu}

We presented a comparison between the performance of Krylov subspace iterative solvers based on  Multiple Trace Formulations (MTF) and Domain Decomposition Methods (DDM) for two dimensional frequency domain scattering problems involving bounded composite materials,i.e.~piece-wise constant material parameters. We investigated DDM based both on classical and generalized/DtN Robin boundary conditions. The former DDM version gives rise to linear systems which, although not particularly well-suited for iterative solvers, can be solved efficiently via Schur complements. On the other hand, the generalized DDM version that incorporates square root approximations of DtN maps gives rise to relatively small numbers of iterations. Extensions to three-dimensional configurations are currently underway. 

\section*{Acknowledgments}
Catalin Turc gratefully acknowledges support from NSF through contracts DMS-1312169 and DMS-1614270. Carlos Jerez-Hanckes thanks partial support from Conicyt Anillo ACT1417.

\section{Appendix}

 In this section we consider DDM for the Helmholtz equation in one dimension. More precisley, we consider the Helmholtz equation
 \begin{eqnarray}\label{eq:1D}
   u''(x)+(k(x))^2u(x)&=&0\quad {\rm in}\quad (a,b)\nonumber\\
   u(a) = A\quad &{\rm and}& u(b)=B
 \end{eqnarray}
 where the wavenumber $k(x)$ is a piecewise constant function, that is
 \[
 k(x)=k_j\quad x\in (a_j,a_{j+1}),\quad \cup_{j=0}^{N+1}[a_j,a_{j+1}]=[a,b]
 \]
 and $u$ and $u'$ are continuous at $a_j,j=0,\ldots,N+1$. We note that we do not require that the wavenumbers $k_j$ be necessarily different on adjacent intervals. The classical DDM formulation of the Helmholtz equation above can be written in the form
 \begin{eqnarray}
   u_j''+k_j^2 u_j&=&0 \quad{\rm in}\quad (a_j,a_{j+1})\nonumber\\
   f_{j,j-1}:=(-u_j'+i\eta\ u_j)|_{x=a_j}&=&(-u_{j-1}'+i\eta\ u_{j-1})|_{x=a_j}\nonumber\\
   f_{j,j+1}:=(u_j'+i\eta\ u_j)|_{x=a_{j+1}}&=&(u_{j+1}'+i\eta\ u_{j+1})|_{x=a_{j+1}}\nonumber
 \end{eqnarray}
 for all $1\leq j\leq N$ together with the end-interval equations
 \begin{eqnarray}
   u_0''+k_0^2 u_0&=&0 \quad{\rm in}\quad (a_0,a_{1})\nonumber\\
   u_0(a_0)&=& A\nonumber\\
   f_{0,1}:=(u_0'+i\eta\ u_0)|_{x=a_{1}}&=&(u_{1}'+i\eta\ u_{1})|_{x=a_{1}}\nonumber
 \end{eqnarray}
 and
\begin{eqnarray}
   u_{N+1}''+k_{N+1}^2 u_{N+1}&=&0 \quad{\rm in}\quad (a_{N+1},a_{N+2})\nonumber\\
   f_{N+1,N}:=(-u_{N+1}'+i\eta\ u_{N+1})|_{x=a_{N+1}}&=&(-u_{N}'+i\eta\ u_{N})|_{x=a_{N+1}}\nonumber\\
  u_{N+1}(a_{N+2})&=&B.\nonumber
\end{eqnarray}
To each of these Robin problems we associate RtR maps. First we define for $1\leq j\leq N$ the following matrices
\begin{equation}
  \mathcal{S}^j\begin{bmatrix}f_{j,j-1}\\f_{j,j+1}\end{bmatrix}:=\begin{bmatrix}(u_j'+i\eta\ u_j)|_{x=a_j}\\(-u_j'+i\eta\ u_j)|_{x=a_{j+1}}\end{bmatrix}\nonumber
\end{equation}
then the following complex scalars
\begin{equation}
  \mathcal{S}^0f_{0,1}=(-u_0'+i\eta\ u_0)|_{x=a_{1}}+\gamma_0 A\nonumber
\end{equation}
and
\begin{equation}
  \mathcal{S}^{N+1}f_{N+1,N}=(u_{N+1}'+i\eta\ u_{N+1})|_{x=a_{N+1}}+\gamma_{N+1}B\nonumber.
  \end{equation}
Denoting $h_j=a_{j+1}-a_j$, it is a straightforward matter to compute
\[
\mathcal{S}^j_{11}=\mathcal{S}^j_{22}=\frac{(k_j+\eta)^2(e^{ik_jh_j}-e^{-ik_jh_j})}{(k_j-\eta)^2e^{-ik_jh_j}-(k_j+\eta)^2e^{ik_jh_j}}
\]
and
\[
\mathcal{S}^j_{12}=\mathcal{S}^j_{21}=-\frac{4k_j\eta}{(k_j-\eta)^2e^{-ik_jh_j}-(k_j+\eta)^2e^{ik_jh_j}}
\]
for $1\leq j\leq N$. We also get
\[
\mathcal{S}^0=\frac{(\eta+k_0)e^{-ik_0h_0}-(k_0-\eta)e^{ik_0h_0}}{(\eta-k_0)e^{-ik_0h_0}-(k_0+\eta)e^{ik_0h_0}}
\]
and
\[
\mathcal{S}^{N+1}=\frac{(\eta+k_{N+1})e^{-ik_{N+1}h_{N+1}}+(k_{N+1}-\eta)e^{ik_{N+1}h_{N+1}}}{(\eta-k_{N+1})e^{-ik_{N+1}h_{N+1}}-(k_{N+1}+\eta)e^{ik_{N+1}h_{N+1}}}.
\]
Ordering the data $f=[f_{01}\ f_{10}\ f_{12}\ \ldots f_{N+1,N}]^\top$, then the classical DDM can be written in the form $(I+A)f=g$ where the matrix $I+A$ is given in explicit form
\begin{equation}\label{Aexp1D}
  I+A=\begin{bmatrix}I&-\mathcal{S}^1_{11}& -\mathcal{S}^1_{12}&0&0&0&0&0&\ldots&0\\
  -\mathcal{S}^0& I& 0& 0 & 0 & 0 & 0& 0 & \ldots& 0\\
  0 & 0 & I & -\mathcal{S}^2_{11} & -\mathcal{S}^2_{12} & 0 & 0 & 0 & \ldots &0\\
  0 & -\mathcal{S}^1_{21} & -\mathcal{S}^1_{22} & I & 0 & 0 & 0 & 0 & \ldots & 0\\
  0 & 0 & 0 & 0 & I & -\mathcal{S}^3_{11} & -\mathcal{S}^3_{12} & 0 & \ldots & 0\\
  \ldots & \ldots & \ldots & \ldots & \ldots & \ldots & \ldots & \ldots &\ldots &\ldots\\
  0 & 0 & 0 & 0 & 0 & 0& 0 & \ldots  & I & -\mathcal{S}^{N+1}\\
  0 & 0 & 0 & 0 & 0 & 0 & \ldots & -\mathcal{S}^N_{21} & -\mathcal{S}^N_{22} & I
  \end{bmatrix}.
\end{equation}

\begin{figure}
\centering
\includegraphics[scale=0.32]{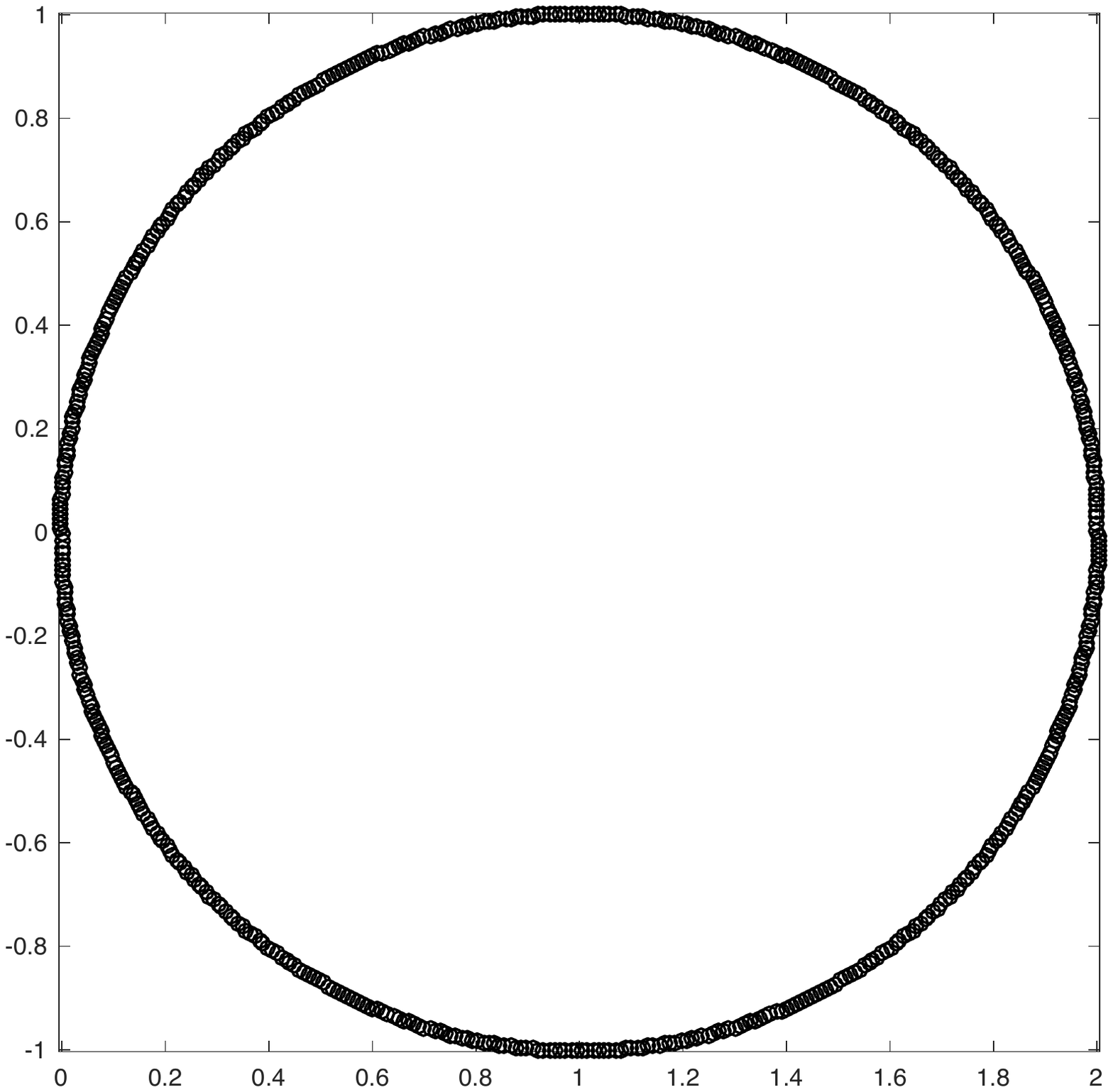}
\includegraphics[scale=0.4]{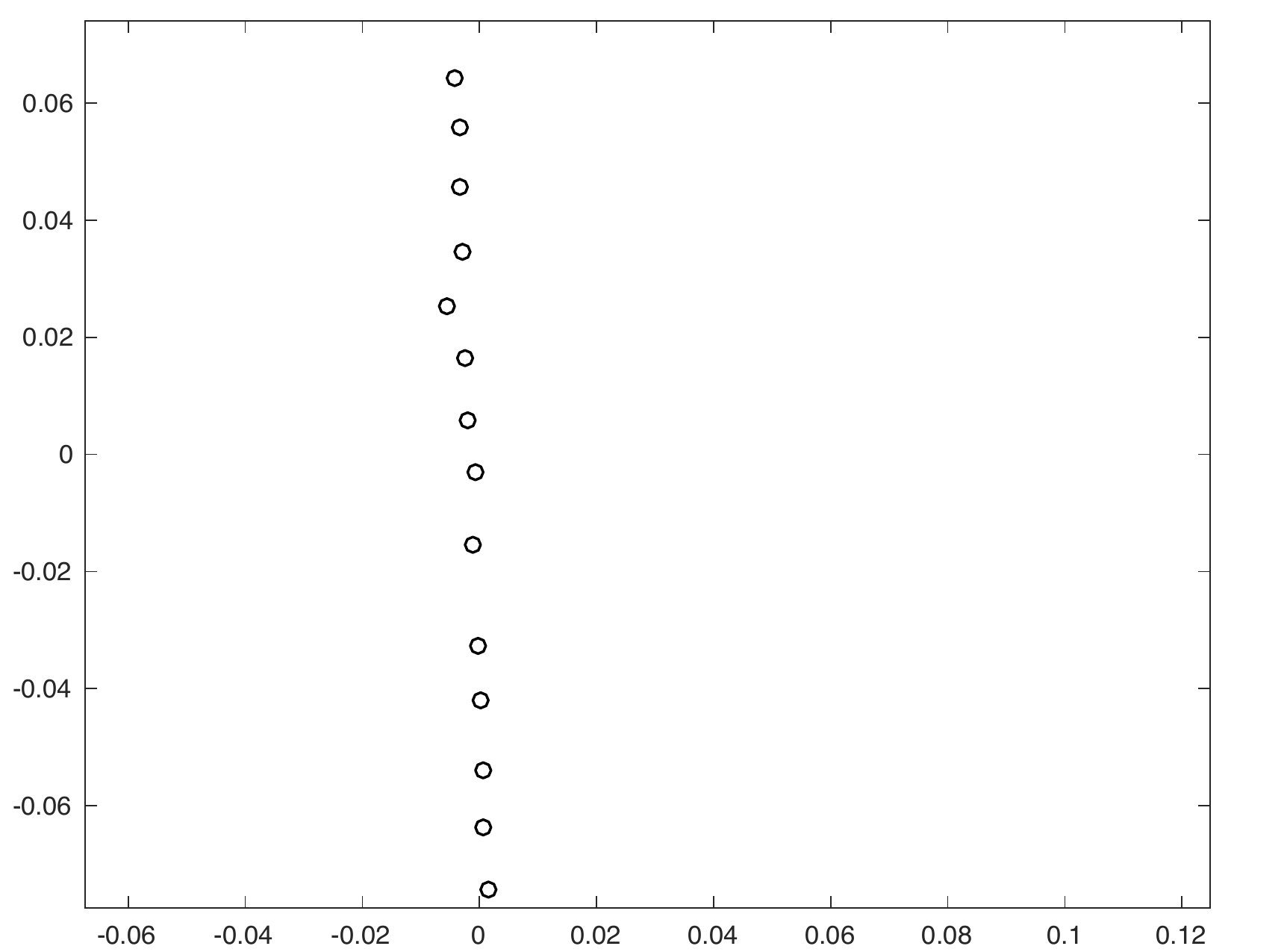}
\caption{Distribution of eigenvalues of the matrix $I+A$ defined in equation~\eqref{Aexp1D} for $\eta=1$ in the case when we solve the Helmholtz equation on the interval $[0,1]$ with $k_0=1$ in $(0,1/4)$, $k_1=2$ in $(1/4/,1/2)$, $k_2=4$ in $(1/2,3/4)$, and $k_3=8$ in $(3/4,1)$. We further subdivided the interval $(0,1/4)$ into 20 subintervals of equal length, the interval $(1/4/,1/2)$ into 40 subintervals of equal length, the interval $(1/2,3/4)$ into 80 subintervals of equal length, and finally the interval $(3/4,1)$ into 160 subintervals of equal length. The smallest eigenvalues of the ensuing matrix $I+A$ is of the order $10^{-3}$.}
\label{fig:eig}
\end{figure}

We present in Figure~\ref{fig:eig} the spectral properties of the matrix $I+A$ defined in equation~\eqref{Aexp1D} for a case of piecewise constant wavenumber that takes four values in the interval $(0,1)$, and a total of 300 subintervals.  The spectral properties of the ensuing DDM are associated with poor behavior of GMRES iterative solvers: the eigenvalues are distributed almost uniformly on a circle of radius close to one centered at $(1,0)$.  

In the case of DtN DDM algorithm we make use of the following DtN maps, assumed to be properly defined:
\[
-v'_j(a_j)=dtn^{-}(a_j) v_j(a_j)
\]
where $v_j$ is the solution of the following problem
\begin{eqnarray}
   v_j''+k_j^2 v_j&=&0 \quad{\rm in}\quad (a_j,a_{j+1})\nonumber\\
   v_j(a_j)=A_j,&\quad& v_j(a_{j+1})=0\nonumber\\
\end{eqnarray}
and
\[
w'_j(a_{j+1})=dtn^{+}(a_{j+1}) w_j(a_{j+1})
\]
where $w_j$ is the solution of the following problem
\begin{eqnarray}
   w_j''+k_j^2 w_j&=&0 \quad{\rm in}\quad (a_j,a_{j+1})\nonumber\\
   w_j(a_j)=0,&\quad& w_j(a_{j+1})=A_{j+1}.\nonumber
\end{eqnarray}
It can be easily shown that
\[
dtn^{-}(a_j)=dtn^{+}(a_{j+1})=-ik_j\frac{e^{ik_jh_j}+e^{-k_jh_j}}{e^{-ik_jh_j}-e^{ik_jh_j}}.
\]
The DtN DDM formulation of the Helmholtz equation above can be written in the form
 \begin{eqnarray}
   u_j''+k_j^2 u_j&=&0 \quad{\rm in}\quad (a_j,a_{j+1})\nonumber\\
   f^{dtn}_{j,j-1}:=(-u_j'+dtn^{+}(a_j) u_j)|_{x=a_j}&=&(-u_{j-1}'+dtn^{+}(a_j)\ u_{j-1})|_{x=a_j}\nonumber\\
   f^{dtn}_{j,j+1}:=(u_j'+dtn^{-}(a_{j+1})\ u_j)|_{x=a_{j+1}}&=&(u_{j+1}'+dtn^{-}(a_{j+1})\ u_{j+1})|_{x=a_{j+1}}\nonumber
 \end{eqnarray}
 for all $1\leq j\leq N$ together with corresponding end-interval equations. Corresponding RtR DtN maps/matrices can be defined and their entries are given by
 \[
 \mathcal{S}^{dtn,j}_{11}=\mathcal{S}^{dtn,j}_{22}=\frac{4k_j^2}{d_j(e^{-ik_jh_j}-e^{ik_jh_j})}.
 \]
 and
 \begin{eqnarray*}
   \mathcal{S}^{dtn,j}_{12}&=&-2ik_j\frac{dtn^{-}(a_j)+dtn^{+}(a_j)}{d_j}\\
   \mathcal{S}^{dtn,j}_{21}&=&-2ik_j\frac{dtn^{+}(a_{j+1})+dtn^{-}(a_{j+1})}{d_j}
 \end{eqnarray*}
 where
 \begin{equation*}
   d_j=(dtn^{-}(a_j)-ik_j)(dtn^{+}(a_{j+1})-ik_j)e^{-ik_jh_j}-(dtn^{-}(a_j)+ik_j)(dtn^{+}(a_{j+1})+ik_j)e^{ik_jh_j}.
 \end{equation*}
With these notation the DtN DDM can be written in the form $(I+A^{dtn})f^{dtn}=g^{dtn}$ where the matrix $I+A^{dtn}$ is given in explicit form
\begin{equation}\label{Aexp1Ddtnf}
  I+A^{dtn}=\begin{bmatrix}I&-\mathcal{S}^{dtn,1}_{11}& -\mathcal{S}^{dtn,1}_{12}&0&0&0&0&0&\ldots&0\\
  0 & I& 0& 0 & 0 & 0 & 0& 0 & \ldots& 0\\
  0 & 0 & I & -\mathcal{S}^{dtn,2}_{11} & -\mathcal{S}^{dtn,2}_{12} & 0 & 0 & 0 & \ldots &0\\
  0 & -\mathcal{S}^{dtn,1}_{21} & -\mathcal{S}^{dtn,1}_{22} & I & 0 & 0 & 0 & 0 & \ldots & 0\\
  0 & 0 & 0 & 0 & I & -\mathcal{S}^{dtn,3}_{11} & -\mathcal{S}^{dtn,3}_{12} & 0 & \ldots & 0\\
  \ldots & \ldots & \ldots & \ldots & \ldots & \ldots & \ldots & \ldots &\ldots &\ldots\\
  0 & 0 & 0 & 0 & 0 & 0& 0 & \ldots  & I & 0\\
  0 & 0 & 0 & 0 & 0 & 0 & \ldots & -\mathcal{S}^{dtn,N}_{21} & -\mathcal{S}^{dtn,N}_{22} & I
  \end{bmatrix}.
\end{equation}
 
In what follows we present a strategy to construct an effective preconditioner for the linear systems associated with the classical DDM formulations. The main idea is to replace the entries in the matrix $I+A^{dtn}$ corresponding to $\mathcal{S}^{dtn,j}_{11}$ and $\mathcal{S}^{dtn,j}_{22}$ by zero. At this stage we find more intuitive to refer to $\mathcal{S}^{dtn,j}_{12}$ as to $\mathcal{S}^{dtn,j}_b$ (the subscript stands for backward, consistent with the direction in which the information propagates) and to $\mathcal{S}^{dtn,j}_{21}$ as to $\mathcal{S}^{dtn,j}_f$ (the subscript stands for forward). We consider thus the following matrix
\begin{equation}\label{Aexp1Dtn}
  I+\widetilde{A}^{dtn}=\begin{bmatrix}I&0& -\mathcal{S}^{dtn,1}_b&0&0&0&0&0&\ldots&0\\
  0& I& 0& 0 & 0 & 0 & 0& 0 & \ldots& 0\\
  0 & 0 & I & 0 & -\mathcal{S}^{dtn,2}_b & 0 & 0 & 0 & \ldots &0\\
  0 & -\mathcal{S}^{dtn,1}_f & 0 & I & 0 & 0 & 0 & 0 & \ldots & 0\\
  0 & 0 & 0 & 0 & I & 0 & -\mathcal{S}^{dtn,3}_b & 0 & \ldots & 0\\
  \ldots & \ldots & \ldots & \ldots & \ldots & \ldots & \ldots & \ldots &\ldots &\ldots\\
  0 & 0 & 0 & 0 & 0 & 0& 0 & \ldots  & I & 0\\
  0 & 0 & 0 & 0 & 0 & 0 & \ldots & -\mathcal{S}^{dtn,N}_f & 0 & I
  \end{bmatrix}.
\end{equation}
The matrices $I+\widetilde{A}^{dtn}$ corresponding to the same experiment described in Figure~\ref{fig:eig} have only one eigenvalue $\lambda=1$ with algebraic multiplicity $2(N+1)$ (this is the number of unknown in the DDM) and geometric multiplicity 2, that is it has only two linearly independent eigenvectors, which turn out to be the first and the last canonical vectors in $\mathbb{R}^{2(N+1)}$. This situation was already pointed out in~\cite{vion2014double} in the case of constant wavenumber. The key property that facillitates the use of the matrix $I+\widetilde{A}^{dtn}$ as the basis of a DDM preconditioner is the fact that the inverse of this matrix can be computed explicitly, and its expression does not involve algebraic inverses.  
\begin{figure}
\centering
\includegraphics[height=50mm]{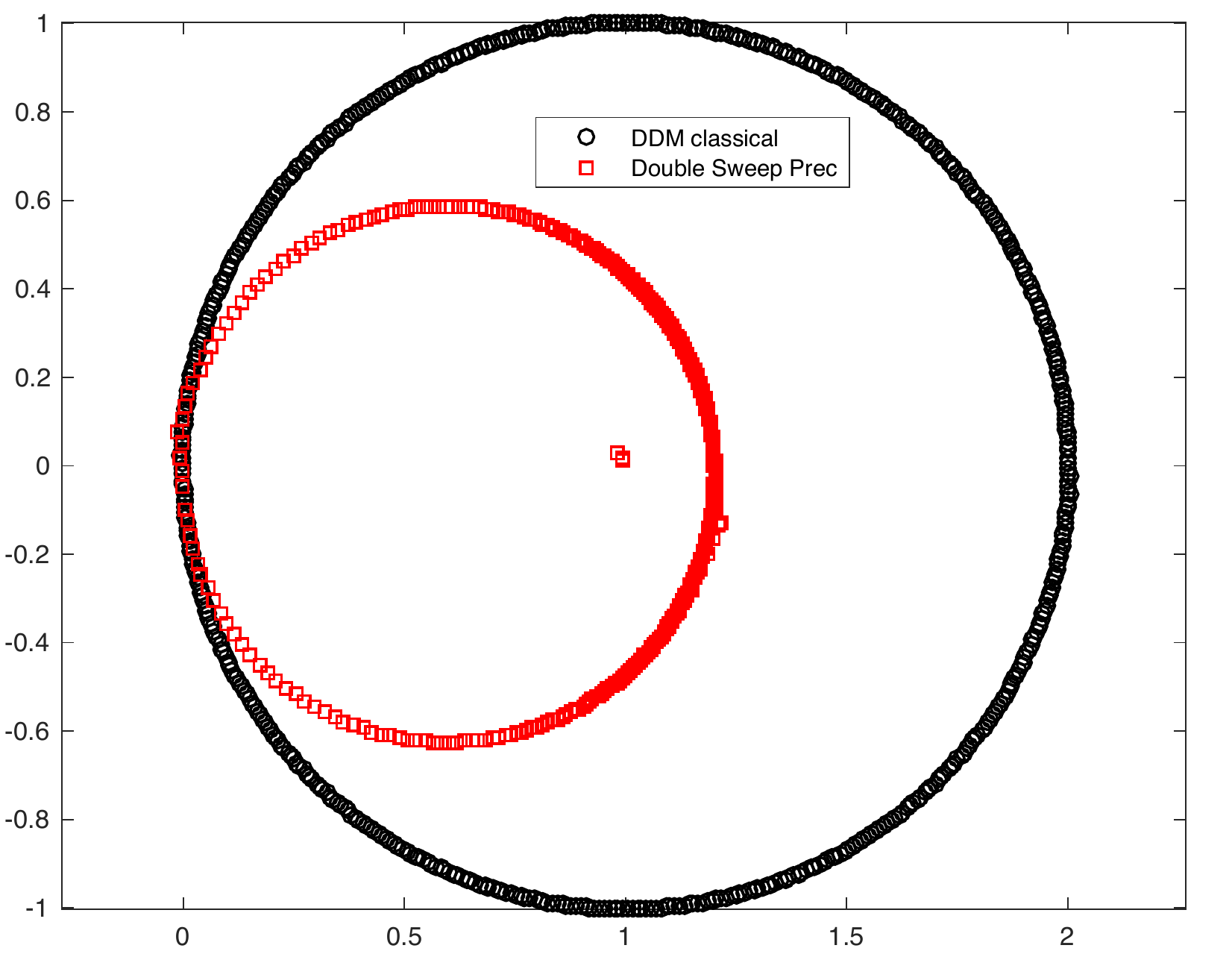}
\caption{Distributions of eigenvalues of the matrices $I+A$, $I+A^{dtn}$, and $(I+\widetilde{A}^{dtn})^{-1}(I+A)$ respectively.}
\label{fig:eig2}
\end{figure}

\bibliography{biblioTJ}

\end{document}